\documentclass{amsart}
\usepackage{graphicx}
\usepackage{latexsym,amscd,amssymb}
\usepackage{amsfonts}

\usepackage{graphicx}

\theoremstyle{plain}
   \newtheorem{theorem}{Theorem}[section]
   \newtheorem{proposition}[theorem]{Proposition}
   \newtheorem{lemma}[theorem]{Lemma}
   \newtheorem{corollary}[theorem]{Corollary}
   
\theoremstyle{definition}
   \newtheorem{definition}[theorem]{Definition}
   \newtheorem{example}[theorem]{Example}

   \newtheorem{remark}[theorem]{Remark}

\numberwithin{equation}{section}

\newcommand\rp{{\overset{p_r}{\sim}}}
\newcommand\rcp{{\overset{p_r^*}{\sim}}}

\newcommand\sh{{\mathrm{shape}}}
\newcommand\lb{{\mathrm{label}}}
\newcommand\rda{{\overset{ \mathrm{D}_1^r}{\sim}}}
\newcommand\rdb{{\overset{ \mathrm{D}_2^r}{\sim}}}
\newcommand\rdc{{\overset{ \mathrm{D}_3^r}{\sim}}}

\newcommand\Knuth{{\overset{K}{\sim}}}
\newcommand\dKnuth{{\overset{K^*}{\sim}}}

\newcommand\Des{{\mathrm{Des}}}

\newcommand\al{\alpha}

\newcommand{\cellsize}{16}
\newlength{\cellsz} 
\newsavebox{\cell}
\sbox{\cell}{\begin{picture}(\cellsize,\cellsize)
\put(0,0){\line(1,0){\cellsize}}
\put(0,0){\line(0,1){\cellsize}}
\put(\cellsize,0){\line(0,1){\cellsize}}
\put(0,\cellsize){\line(1,0){\cellsize}}
\end{picture}}
\newcommand\cellify[1]{\def\thearg{#1}\def\nothing{}%
\ifx\thearg\nothing
\vrule width0pt height\cellsz depth0pt\else
\hbox to 0pt{\usebox{\cell} \hss}\fi%
\vbox to \cellsz{
\vss
\hbox to \cellsz{\hss$#1$\hss}
\vss}}
\newcommand\tableau[1]{\vtop{\let\\\cr
\baselineskip -16000pt \lineskiplimit 16000pt \lineskip 0pt
\ialign{&\cellify{##}\cr#1\crcr}}}

\begin{document}

\title[plactic relations for $r$-domino tableaux]
{Plactic relations for $r$-domino tableaux}
\author{ M\"{U}GE TA\c{S}KIN}
\address{BO\~{G}AZ\.{I}\c{C}\.{I} \"{U}N\.{I}VERS\.{I}TES\.{I}}
\thanks{This research was partially supported by the Fields Institute and York University, Toronto, ON, CA.}

\begin{abstract}
The  work of C. Bonnaf{\'e}, M.Geck, L. Iancu and T. Lam
\cite{Geck-Lam} shows through two conjectures that $r$-domino
tableaux have an important role   in Kazhdan-Lusztig theory  of type
$B$ with unequal parameters. In this paper we provide plactic
relations on signed permutations which determine whether given two
signed permutations have the same insertion $r$-domino tableaux in
Garfinkle's algorithm \cite{Garfinkle1}. Moreover, we show that a
particular extension of these relations can describe Garfinkle's
equivalence relation \cite{Garfinkle1} on $r$-domino tableaux which
is given through the notion of  open cycles. With these results we
enunciate the conjectures of  \cite{Geck-Lam} and provide necessary
tool for their proofs.
\end{abstract}

\maketitle


\section{Introduction} \label{Introduction}

Let $W$ be a finite Coxeter  group and let  $L:W\mapsto
\mathbb{Z}_{\geq 0}$ be a weight function such that
$$L(uw)=L(u)+L(w) ~\text{ if and only if }~ l(uw)=l(u)+l(w)$$ where
$l:W\mapsto  \mathbb{Z}_{\geq 0}$ is the usual length function on
$W$.  As it is described by Lusztig in \cite{Lusztig} every weight
function determines an Iwahori-Hecke algebra  and  three preorders
on $W$ whose equivalence classes are called {\it left, right} and
{\it two-sided cells}.  The  importance of these cells lies in the
fact that  they carry representations of $W$ and its corresponding
Iwahori-Hecke algebra $\mathcal{H}$. Furthermore they have an
important role in the representation theory of reductive algebraic
groups over finite or $p$-adic fields \cite{Lusztig} and in the
study of rational Cherednik algebras \cite{Gordon} and the
Calogero-Moser spaces \cite{Gordon2}.

The case $L=l$ is in fact first introduced by Kazhdan and Lusztig in \cite{Kazhdan-Lusztig}
as a purely combinatorial
 tool for the theory of primitive ideals in the universal enveloping algebras of semisimple
 complex Lie algebras. In
 this case  the combinatorial characterizations of  cells are   well known, where Knuth
 (or plactic)  relations appear as the mediating tool. Namely,  when $W$ is  type $A$ then
  each  right (left) cell  corresponds to the plactic (respectively coplactic)  class of some standard
Young tableau, whereas each two-sided cell consists of those
permutations  which lie in the plactic classes of tableaux of the
same shape.  This characterizations depend on  Joseph's
classification  of  primitive ideals  in type A,  where  Knuth
(plactic) relations play a crucial role.

In the types B, C and D, on the other hand the emerging
combinatorial objects are standard domino tableaux. The connection
is first revealed in the work of Barbash and Vogan
\cite{Barbash-Vogan} where they provide necessary conditions for the
characterizations of primitive ideals  through   an algorithm which
uses the palindrome representations of signed permutations in order
to assign to  every signed permutation $\alpha$  a pair of same
shape standard $r$-domino tableaux $(P^r(\alpha), Q^r(\alpha))$
bijectively, for $r=0$ or $r=1$. Meanwhile,  an analog of Knuth
relations provided by Joseph in \cite{Joseph} established the
sufficient conditions. On the other hand  Garfinkle
\cite{Garfinkle1,Garfinkle2,Garfinkle3} finalized the classification
problem for these types  by showing through her two algorithms on
domino tableaux that these two sets of relations are in fact
equivalent.  Her first algorithm   assigns any signed permutation to
a pair of same shape standard  $r$-domino tableaux for $r$ equal to
$0$ or $1$ and the second defines an equivalence relation between
domino tableaux through the notion of {\it open cycles}.  We remark
that the extension of Garfinkle and Barbash-Vogan algorithm  for
larger $r$ is given in \cite{Leeuwen} and \cite{Geck-Lam}
respectively.

The case $L\not = l$ is also known as unequal parameter
Kazhdan-Lusztig theory and it appears for the types $B_n$, $I_2(n)$
and $F_4$, where the classification problem for the latter two can
be dealt with computational methods, see \cite{Geck}. For type
$B_n$, the weight  function is determined by two integers $a,b>0$
such that
$$L(s_i)=\left \{ \begin{aligned} a  &~\text{ if }~ 1\leq i\leq n-1 \\ b &~\text{ if }~ i=0
\end{aligned} \right.
$$
where $s_0$ is the transposition  $(-1,1)$ and $\{s_i=(i,i+1) |
1\leq i\leq n-1\}$ are the type $A$  generators of $B_n$. Recently,
the role of $r$-domino tableaux  in this theory is revealed in  the
work of Bonnaf{\'e}, Geck, Iancu, and Lam \cite{Geck-Lam} through
two main conjectures:
\begin{enumerate}
\item[ $\bullet$] {\it Conjecture A:}
 If  $ r a < b < (r+1)a$ for some $r\geq 0$ then  two signed permutations lie in  the same
 Kazhdan-Lusztig right (left) cell
 if and only if their  insertion (recording) $r$-domino tableau are the same.
\item[ $\bullet$]{\it Conjecture B:}  If   $b=r a$  for some $r \geq 1$ then
two signed permutations lie in the same  Kazhdan-Lusztig right
(left) cell if and only if their  insertion (recording) $r-1$-domino
tableau are equivalent through the notion of open cycles.
\end{enumerate}

In order to establish the proofs of these conjecture one definitely
needs the plactic relations between signed permutations which
determines when the insertion $r$-domino tableaux of two signed
permutations are the same  or equivalent through the notion of open
cycles. Our aim here is to fill this gap.

This paper is organized as follows:  The descriptions of
Barbash-Vogan and Garfinkle's algorithms can be found in
Section~\ref{related.background} together with some lemmas which are
essential in the following section. In
Section~\ref{proof.of.main.theo.sec} the definition  of plactic
relations are given and  they are shown to be necessary and
sufficient  for describing plactic classes of $r$-domino tableaux.

\begin{remark}  Recently T. Pietraho \cite{Pietraho3}  has
found another set of generators which can be shown to be equivalent to
$\mathrm{D}_1^r$, $\mathrm{D}_2^r$, $\mathrm{D}_3^r$ and
$\mathrm{D}_3^{r-1}$ given in the
Definition~\ref{plactic.relations}. On the other hand these
relations describes a larger set, namely the set of all permutations
whose insertion $r$-domino tableaux are equivalent through the
notion of open cycles. Finally, by using his results and an earlier
version of the present work, C. Bonnaf\'{e} provides a partial
result towards the previous conjectures \cite{Bonnafe}.
\end{remark}

\section*{Acknowledgments} The author grateful to referee for his/her crucial comments and  to Nantel
 Bergeron, Victor Reiner and Huilan Li for helpful discussions and comments.

\section{Related background} \label{related.background}

A sequence $\lambda=(\lambda_1, \ldots, \lambda_k)$  is a {\it
partition} of $n$, denoted by  $\lambda \vdash n$,  if
$\sum_{i=1}^{k} \lambda_i=n ~\text{ and }~ \lambda_{i}\geq
\lambda_{i+1} >0$ where its   {\it Ferrers diagram} consists of left
justified arrows of boxes such that the $i$-th row has $\lambda_i$
boxes. For example
$$ \lambda =(2,2,1)= ~\tableau{{ \
}&{ \ }\\ { \ }&{ \ } \\{ \ }}  $$ A partition $\lambda=(\lambda_1,
\ldots, \lambda_k)$ can be also seen as a set of integer pairs
$(i,j)$ such that $1\leq i \leq k$ and $1\leq j \leq \lambda_i$.
Therefore for two partitions    $\lambda$ and $\mu$,  we can  define
usual set operations such as $\lambda\cup \mu$, $\lambda\cap \mu$,
$\lambda\subset \mu$, $\lambda-\mu$,  but the  resulting sets do not
necessarily  correspond to any partitions.

\begin{definition} \label{partition.def.1}
For  two partitions  $\lambda$ and $\mu$ satisfying $\mu\subset \lambda$ we define
$ \lambda/\mu= \lambda-\mu$
to be  the {\it skew partition} determined by  $\lambda$ and $\mu$.
\end{definition}

\begin{definition} \label{partition.def.2}
 Let $\gamma$ and $\gamma'$ be two skew shapes.
\begin{enumerate}
\item[1.] If  $\gamma\cap \gamma'=\emptyset$ and $\gamma \cup
\gamma'$ also  corresponds   a skew shape  then we define
$\gamma \oplus\gamma'=\gamma \cup \gamma'$.
\item[2.] If  $\gamma' \subset \gamma$ and $\gamma -
\gamma'$ also  corresponds   a skew shape  then we define
$\gamma \ominus  \gamma'=\gamma- \gamma'$.
\end{enumerate}
\end{definition}

\begin{definition}  \label{partition.def.3}
Let   $\lambda$ be   a partition and  $(i,j) \in \lambda$.
\begin{enumerate}
\item[1.] If $ (i,j) \in  \lambda$ and $ \lambda \ominus {(i,j)}$ is also a partition  then
 $(i,j)$ is called a {\it corner }of $\lambda$.
 \item[2.] If $ (i,j)\not \in  \lambda$ and $
\lambda \oplus {(i,j)}$ is also a partition  then $(i,j)$ is called
an {\it empty corner} of $\lambda$.
\end{enumerate}
\end{definition}

\begin{definition}  \label{tableau.def.1}
 A {\it  skew tableau} $T$  of shape $\lambda/\mu$ is obtained by labeling the cells  of  $\lambda/\mu$
with non repeating, totally ordered letters such that the letters increase
from left to right and from top to bottom.
 If $\mu =\emptyset$ then  $T$ is called a {\it Young tableau}. We denote by
$$\lb(T) ~\text{ and }~\sh(T)
$$
respectively, the set of letters labeling  each box of $T$ and  the
partition underlying $T$. For a set $\mathcal{A}$ of letters given
with a totally ordering, we denote by $SYT_{\mathcal A}$ the set of
Young tableaux labeled with $\mathcal{A}$. If
$\mathcal{A}=\{1,2,\ldots,n\}$   then $T$ is called a   {\it
standard skew} or {\it standard Young tableau} according to the
shape of $T$. Moreover  the set of standard Young tableaux are
denoted by $SYT_n$.
 \end{definition}

Let $S_{\mathcal A}$ ($S_n$) denote  the symmetric group on the
totaly ordered set $\mathcal{A}$ (respectively on $\{1,2,\ldots,
n\}$). When the size of $\mathcal A$ is $n$ we have an order
preserving bijection between $\mathcal{A}$ and $\{1,2,\ldots, n\}$
and this yields two more  bijections between $S_{\mathcal A}$ and
$S_n$ as well as $SYT_{\mathcal A}$ and $SYT_n$. Therefore the
following discussions  and results can be generalized to any finite
totally  ordered set $\mathcal{A}$.

There is an important connection,  between  Young tableaux $SYT_n$
and the symmetric group $S_n$, known as the  {\it Robinson-Schensted
correspondence (RSK)}, which was realized by Robinson and Schensted
independently. In this correspondence, every permutation $w \in
S_{n}$ is assigned bijectively  to a pair of same shape tableaux
$(P(w),Q(w))$ in $SYT_n\times SYT_n$ through {\it insertion} and
{\it recording} algorithms. Let us explain these algorithms briefly.
We denote by $(P_{i-1},Q_{i-1})$ the  tableaux obtained by insertion
and recording algorithms on the first $i-1$ indices of $w=w_1\ldots
w_n$. In order to get $P_i$ we proceed as follows: if $w_{i}$ is
greater then the last number on the first row of $P_{i-1}$, then
$w_i$ is concatenated to the first row of $P_{i-1}$ from the right,
otherwise $w_i$ replaces the smallest number, say $a$, among all
numbers in the first row which are greater then $w_i$ and the
insertion algorithm continues with the insertion of  $a$ to  the
next row. Observe that after finitely many steps the insertion
algorithm terminates with a new appearing cell on some row of
$P_{i-1}$. The resulting tableau is then $P_i$ and the recording
tableau $Q_i$ is found by filling this new cell in $Q_{i-1}$ with
the number $i$. We illustrate these algorithms with the following
example.

\begin{example}\label{tableau.ex.1}
 Let $w=52413 \in S_5$. Then,
$$\begin{aligned}
P_1=5, ~&~ P_2=\begin{array}{c}2\\5\end{array}, ~&~
P_3=\begin{array}{cc}2&4\\5\end{array}, ~&~
P_4=\begin{array}{cc}1&4\\2\\5\end{array}, ~&~
P_5=\begin{array}{cc}1&3\\2&4\\5\end{array}=P(w)&\\
Q_1=1, ~&~ Q_2=\begin{array}{c}1\\2\end{array}, ~&~
Q_3=\begin{array}{cc}1&3\\2\end{array}, ~&~
Q_4=\begin{array}{cc}1&3\\2\\4\end{array}, ~&~
Q_5=\begin{array}{cc}1&3\\2&5\\4\end{array}=Q(w)&
\end{aligned}
$$
\end{example}

The following result of Sch\"utzenberger \cite{Schutzenberger2}
reveals an important property of the RSK.

\begin{theorem}  \label{tableau.theo.0} If $w\in S_n$, then
$$P(w^{-1})=Q(w) ~\text{and}~ Q(w^{-1})=P(w).$$
\end{theorem}

 There are two equivalence relations and a related result
given  by Knuth \cite{Knuth} which are  fundamental  in the
combinatorics of tableaux. In the following we provide them in a
more general setting:

\begin{definition} \label{tableau.def.2} Let $\mathcal{A}$ be a
totally ordered set of letters and  $u=u_1 \ldots u_n \in
S_\mathcal{A}$. If either  $ u_i <u_{i+2}<u_{i+1}$ or  $ u_{i}
<u_{i-1}< u_{i+1}$ for some $i$ then
$$u=u_1 \ldots u_{i-1} (u_{i} ~u_{i+1}) ~u_{i+2} \ldots u_n~ \sim~ u_1 \ldots u_{i-1} (u_{i+1} ~u_{i}) ~u_{i+2} \ldots u_n=u'. $$

We say  $u, w \in S_{\mathcal{A}}$ are {\it Knuth equivalent},  $u
~\Knuth~ w$, if $w$ can be obtained from $u$ by applying a sequence
of $\sim$ relations.   On the other hand  if $u^{-1} ~\Knuth~
w^{-1}$ then $u$ and $w$ are called {\it dual Knuth equivalent}, $u
~\dKnuth~ w$.
\end{definition}

\begin{theorem}[Knuth \cite{Knuth}] \label{tableau.theo.1}
Let $u, w \in S_\mathcal{A}$. Then\\
$i)$ $u  \stackrel{ K}{\cong} w \iff P(u)=P(w)$ \\
$ii)$ $ u \stackrel{K^{*}}\cong w \iff Q(u)=Q( w)$.
\end{theorem}

 We next  illustrate  the  forward and backward
slides of Sch\"utzenberger's {\it jeu de taquin}
\cite{Schutzenberger} without the definition. We remark that  jeu de
taquin slides can be used to give alternative descriptions of both
the Robinson-Schensted algorithm and  Knuth relations. The following
theorem provided by Sch\"utzenberger in \cite{Schutzenberger} reveal
this connection.

\begin{example}\label{tableau.ex.2}
Below we illustrate a forward slide on the tableau $S$ through cell
$c_{12}$ and backward slide on the tableau $T$ through cell
$c_{32}$.
\vskip.1in
$$
S=\tableau{{ \ } & {\bullet}  & {4} \\{ \ }& {2} & {5}  \\{1}& {3}}
\rightarrow \tableau{{ \ }& {2}  & {4} \\
{ \ }&{\bullet}& {5}\\{1}& {3}} \rightarrow \tableau{{ \ }& {2} &
{4}\\{ \ }& {3} & {5}\\{1}& {\bullet}} ~\mbox{ \ \ \ }~ T=\tableau{
{ \ } & {2} & {4} \\ { \ }& {3} & {5} \\{1}& {\bullet}} \rightarrow
\tableau{{ \ } & {2} & {4} \\{ \ }& {\bullet} & {5}
\\{1}& {3}} \rightarrow \tableau{{ \ } & {\bullet}  & {4} \\ { \ }&
{2} & {5}
\\{1}& {3}}$$
\end{example}

\begin{theorem} \label{tableau.theo.2}
 If $P$ is a skew tableau that is brought to a Young
tableau $P'$ by slides, then $P'$ is unique. In fact, $P'$ is the
insertion tableau for the row word  of $P$.
\end{theorem}

\begin{definition}\label{domino.def.6}
The set  of two adjacent cells  $A=\{(i,j),(i,j+1)\}$ (or
$A=\{(i,j),(i+1,j)\}$) is called a horizontal (or respectively
vertical) domino cell.  Now
$$\mathrm{min}(A) ~\mbox{and}~\mathrm{max}(A)$$
denotes the minimum and respectively maximum cell of $A$ in the lexicographic order.
 Let   $A=\{(r_1,c_1),(r_1',c_1')\}$ and  $B=\{(r_2,c_2),(r_2',c_2')\}$  two domino cells in $T$  where $\mathrm{max}(A)=(r_1',c_1')$ and $\mathrm{min}(B)=(r_2,c_2)$.
Then we say
\begin{enumerate}
\item[i)] $B$ lies below $A$  if $\mathrm{min}(B)$ lies below $\mathrm{max}(A)$, equivalently  $r_2>r_1'$.
\item[ii)] $B$ lies to the right of   $A$ if $\mathrm{min}(B)$ lies to the right of  $\mathrm{max}(A)$, equivalently $c_2>c_1'$.
\end{enumerate}

\end{definition}

Let  $\lambda$ be a partition and $A$ be a domino cell.   If
$\lambda \oplus A$ is a partition  then
 $A$ is called an {\it empty
domino corner} of $\lambda$ whereas  if $\lambda \ominus A$ is also
a partition then $A$ is called a {\it  domino corner} of $\lambda$.
Clearly, if a partition has no domino corner then it must be a
$r$-staircase shape $(r,\ldots,2,1)$ for some $r>0$. On the other
hand it is easy to see that  any partition $\lambda$  can be reduced
uniquely  to a $r$-staircase shape $(r,\ldots,2,1)$ for some
$r\geq0$,  by subsequent  removal of  existing domino corners one at
a time.  In this case we say $\lambda$ has  a {\it $2$-core}
equivalent to $(r,\ldots,2,1)$. For $r\geq 0$ we denote by $P(2n,r)$  the set
of all such partitions of size $2n+r(r+1)/2$.

\begin{definition}\label{domino.def.1}
By a {\it labeling} of domino cell $A$ we mean a pair of positive
numbers $(a,a')$ which  label the boxes of $A$ such that   $a\leq
a'$ and $a$ labels  $\mathrm{min}(A)$  and $a'$ labels
$\mathrm{max}(A)$. If the label of $A$ is $(a,a)$ then we say $A$ is
{\it double labeled by} $a$. When we want to indicate the domino
cell $A$ with its labeling,  we use the notation
$$[A,(a,a')]$$
so that $\sh( [A,(a,a')])=A$ and $\lb( [A,(a,a')])=(a,a')$.

 A {\it  $r$-domino tableau} $T$ of shape $\lambda \in P(2n,r)$ is obtained by tiling the
  skew partition  $\lambda/(r,\dots,2,1)$
  with  double labeled horizontal or vertical domino cells $\{[A_1,(a_1,a_1)],
\ldots [A_n,(a_n,a_n)] \}$ such that $a_i>0$ for all $i=1,\ldots,
n$, $a_i \not = a_j$ for $i\not=j$ and  the labels  increase from
left to right and from top to bottom. In this case we have
$$ \lb(T)=\{a_1,a_2,\ldots,a_n\}.$$
A {\it standard $r$-domino tableau} $T$ is a $r$-domino tableau which has  $\lb(T)=\{1,\ldots,n\}$. We denote by  $SDT^r(n)$ the set of all standard $r$-domino tableaux of $n$ dominos.
\end{definition}

\begin{definition}\label{domino.def.2}
 Let $T$ be a $r$-domino tableau and  $\lambda=\sh(T)$. For $A$ is a domino cell in   $\lambda$   and $b\in \lb(T)$ we define,
\begin{enumerate}
\item[1.]  $\mathrm{label}(T, A)$ to be   the pair of integers $(a,a')$ which label the domino cell  $A$ in $T$, where $a\leq a'$.
\item[2.]  $\mathrm{Dom}(T, b)=[B,(b,b)]$  if $B$ is double labeled by $b$ in $T$.
\end{enumerate}
\end{definition}

\begin{example}\label{domino.ex.1}
For example the following is a $2$-domino tableau in
$SDT^2(5)$.\vskip.1in
\[  T = ~\tableau{{ \ }&{ \ }&{1}&{1 }&{5}\\
                   { \ }&{3}&{4}&{4}&{5} \\
                   {2} &{3}\\
                   {2}}  \]
\vskip.1in

Here $T$  has  two domino corners: $A_1=\{(1,5),(2,5)\}$ and
$A_2=\{(2,4),(2,5)\}$, whereas $\mathrm{label}(T, A_1)=(5,5)$ and
$\mathrm{label}(T, A_2)=(4,5)$. On the other hand $\mathrm{Dom}(T,
5)=[A_1,(5,5)]$.
\end{example}


\begin{definition} \label{domino.def.3}
For  two   $r$-domino tableau $S$ and  $T$ satisfying $S\subset T$ we define
$ T/S= T-S$
to be  the {\it skew  $r$-domino tableau} determined by  $S$ and  $T$.
\end{definition}

\begin{definition} \label{domino.def.4}
Let $R$ and $R'$ be two skew $r$-domino tableaux  and let $\sh(R)=\gamma$  and $\sh(R')=\gamma'$.
\begin{enumerate}
\item[1.] If  $\gamma \oplus\gamma'$ is defined and $R\cup R'$ corresponds to a  skew  $r$-domino tableau  as a set  then we define $R\oplus R'=R\cup R'$
\item[2.] If  $\gamma \ominus  \gamma'$ is defined and if $R- R'$ corresponds to a  skew  $r$-domino tableau   as a set then we define $R\ominus R'=R- R'$
\end{enumerate}
\end{definition}

\begin{definition} \label{domino.def.5}
Let $T$ be a (skew) $r$-domino tableau and $a\in \lb(T)$. Then we define
\begin{enumerate}
\item[1.]  $T_{<a}$ ($T_{\leq a}$) to be
the  $r$-domino tableau obtained by restricting $T$ to its double
labeled domino cells whose labels are  less than (and equal to) $a$.
\item[2.] $T_{>a}$ ($T_{\geq a}$) to  be  the skew $r$-domino  tableau obtained
by restricting $T$ to its double labeled domino cells whose labels
are greater than (and equal to) $a$.
\end{enumerate}
\end{definition}


\subsection{ Garfinkle's algorithm}\label{Garfinkle's.algorithm}

Recall that a  signed permutation $\alpha \in B_n$ is a bijection of
$[-n, +n]$ such that $\alpha(-i)=-\alpha(i)$.  The {\it usual
presentation} of  $\alpha \in B_n$ is denoted as
$\alpha=\alpha_1\alpha_2\ldots \alpha_n$ where $\alpha_i=\alpha(i)$
for $1\leq i\leq n$ and $\{|\alpha_1|,|\alpha_2|
\ldots,|\alpha_n|\}=\{1,2,\ldots,n\}$. In the following we set the
following representation for all integers:
$$\bar{a}=\left \{ \begin{aligned} -a &~~\mbox{ if } a>0\\
|a| &~~\mbox{ if } a<0
\end{aligned} \right.
$$


 Garfinkle \cite[Theorem
1.2.13]{Garfinkle1}  provides  an algorithm by which any signed
permutation $\alpha \in B_n$ is assigned bijectively to a pair of
same shape standard $r$-domino tableau $(P^r(\alpha), Q^r(\alpha))$
for $r=0,1$, where $P^r(\alpha)$ is called the {\it insertion} and
$Q^r(\alpha)$ is called the {\it recording} tableau of $\alpha$. Her
algorithm is extended by van Leeuwen \cite{Leeuwen} for  larger
cores.

In the following  we will  explain how to insert an integer into a
$r$-domino tableau according to  Garfinkle's algorithm.  Let $T$ be
a $r$-domino tableau such that $|a| \not \in \lb(T)$.  We denote by
$$T^{\downarrow a}$$
 the tableau which is obtained by inserting $a$ into $T$.

Let $a_0$ be the largest label in $T$  which is smaller then $|a|$.
If $a>0$ then  we first concatenate a horizontal domino labeled with
$(a,a)$ to the first row of $T_{\leq a_0}$ from the right. Otherwise
a vertical domino labeled with $(|a|,|a|)$ is concatenated to the
first column of $T_{\leq a_0}$ from the bottom. Let  $I_0$ denote
the resulting tableau. If the skew tableau $T_{>a_0}$ is empty then
we have
$$T^{\downarrow a}=I_0.$$
Otherwise let $a_1,a_2,\ldots, a_s$ be the increasing sequence of
the labels  in $T_{>a}$. In the following we will find
$T^{\downarrow a}$ through a sequence of tableaux
$I_0,I_1\ldots,I_s$ where
$$T^{\downarrow a}=I_s= I_{s-1}\leftarrow
\mathrm{Dom}(T,a_{s})=\ldots =I_0\leftarrow
\mathrm{Dom}(T,a_{1})\leftarrow \ldots \leftarrow
\mathrm{Dom}(T,a_{s})$$ and for each $i=1,\ldots,s$, $$I_i=
I_{i-1}\leftarrow \mathrm{Dom}(T,a_{i}) $$ is obtained by sliding
$\mathrm{Dom}(T,a_{i})$ to the tableau $I_{i-1}$ in the following
manner: Let
$$B_{i}=\sh(I_{i-1}) \cap \sh(\mathrm{Dom}(T,a_i))
$$
We first assume that
$\mathrm{Dom}(T,a_{i})=[\{(k,l),(k,l+1)\},(a_{i},a_{i})]$ is
horizontal. Then  we have the following possibilities :
\begin{enumerate}
\item[$H_1)$] $B_i=\emptyset$. Then $I_i=I_{i-1}\leftarrow \mathrm{Dom}(T,a_{i})=I_{i-1} \oplus \mathrm{Dom}(T,a_{i})$.
\item[$H_2)$] $ B_i=\{(k,l),(k,l+1)\}$. Then in order to obtain $I_i=I_{i-1}\leftarrow \mathrm{Dom}(T,a_{i})$,  a horizontal domino cell double labeled
by $a_{i}$ is concatenated  to the $(k+1)$-th row of $I_{i-1}$ from
the right.
\item[$H_3)$] $ B_i=\{(k,l)\}$. Then $I_i=I_{i-1}\leftarrow \mathrm{Dom}(T,a_{i})=I_{i-1}\oplus [\{(k,l+1),(k+1,l+1)\},(a_{i},a_i)]$.
\end{enumerate}
\vskip.1in

Now we  assume that
$\mathrm{Dom}(T,a_{i})=[\{(k,l),(k+1,l)\},(a_{i},a_{i})]$ is
vertical. Then we have the following possibilities for $B_i$:
\begin{enumerate}
\item[$V_1)$] $ B_i=\emptyset$. Then $I_i=I_{i-1}\leftarrow \mathrm{Dom}(T,a_{i})=I_{i-1} \oplus \mathrm{Dom}(T,a_{i})$.
\item[$V_2)$] $B_i=\{(k,l),(k+1,l)\}$.  Then in order to find  $I_i=I_{i-1}\leftarrow \mathrm{Dom}(T,a_{i})$, a vertical
domino cell double labeled by $a_{i}$ is concatenated to the
$(l+1)$-th column of $I_{i-1}$ from the bottom.
\item[$V_3)$] $B_i=\{(k,l)\}$. Then $I_i=I_{i-1}\leftarrow \mathrm{Dom}(T,a_{i})=I_{i-1}\oplus [\{(k+1,l),(k+1,l+1)\},(a_{i},a_{i})]$.
\end{enumerate}

Then  insertion and recording $r$-domino tableaux for any
$\alpha=\alpha_1\ldots\alpha_n$ is found in the following way:
Suppose that $P_0$ and $Q_0$ are the tableaux of shape
$(r,\ldots,2,1)$ whose cells are all filed with $0$. For
$\alpha=\alpha_1\ldots\alpha_n \in B_n$ let $P_{i+1}=P_i^{\downarrow
\al_i}$ and let $Q_{i+1}$ be obtained from $Q_i$ by filing the newly
appearing the domino corner of  $P_{i+1}$ with $(i+1,i+1)$ in
$Q_{i+1}$. Then one can obtain $P^r(\al)$ and $Q^r(\al)$ by erasing
all zeros of $P_n$ and respectively $Q_n$.

\vskip.1in The following lemma directly follows from Garfinkle's
algorithm and it indicates some of its  main features.

\begin{lemma} \label{Garfinkle's.algorithm.lemma}Let $n$ be the largest entry in $T$ and  $a$ be an integer satisfying $|a|<n$ and  $a \not \in
\lb(T)$. Then
\begin{enumerate}
\item[$i.$] $T^{\downarrow \overline{a}}=((T^t)^{\downarrow
a})^t$ where $T^t$ is the transpose of $T$.
\item[$ii.$] $T^{\downarrow a}=(T_{<n})^{\downarrow
a}\leftarrow \mathrm{Dom}(T,n)$

\end{enumerate}

\end{lemma}

\begin{example} Find  $T^{\downarrow \overline{2}}$  for
$T=\tableau{ { 1}&{ 1 }&{3}&{3}  \\{4}&{5}&{5}\\{4}&{6}&{6}}$, where
$T_{<2}=\tableau{{1}&{1}}$.

$$\begin{aligned}
&I_0 =\tableau{ {1}&{1} \\ {2} \\{2}}~~
\tableau{ { \ }&{ \ }&{\textbf{3}}&{\textbf{3}}  \\
{4}&{5}&{5}
\\{4}&{6}&{6}}
~\overset{H_1}\longrightarrow~~
I_1=\tableau{ {1}&{1}&{3}&{3} \\ {2}
\\{2}}~~
\tableau{ { \ }&{ \ }&{ \ } \\
{\textbf{4}}&{5}&{5}
\\{\textbf{4}}&{6}&{6}}
~\overset{V_2}\longrightarrow~~
I_3=\tableau{ {1}&{1}&{3}&{3} \\
{2}&{4} \\{2}&{4}}~~
\tableau{ { \ }&{ \ }&{ \ }\\
{ \ }&{\textbf{5}}&{\textbf{5}}
\\{ \ }&{6}&{6}}\\
\\
&\overset{H_3}\longrightarrow~~
I_4=\tableau{ {1}&{1}&{3}&{3} \\
{2}&{4}&{5} \\{2}&{4}&{5}} ~~
\tableau{ { \ }&{ \ }&{ \ } \\
{ \ }&{ \ }&{ \ }
\\{ \ }&{\textbf{6}}&{\textbf{6}}}
~\overset{H_2}\longrightarrow~~ I_5=T^{\downarrow \overline{2}}=\tableau{
{1}&{1}&{3}&{3} \\ {2}&{4}&{5}
\\{2}&{4}&{5}\\{6}&{6}}
\end{aligned}
$$
\end{example}

We now explain the reverse-insertion of domino corners from
$r$-domino tableaux which is the main ingredient of Garfinkle's
bijection. Let $T$ be a $r$-domino tableau and $A$ be a domino
corner in $\sh(T)$. We denote by
$$T^{\uparrow A} ~\text{ and }~\eta(T^{\uparrow A})$$
respectively the  tableau which is  obtained by the
reverse-insertion of $A$, and the number  which is bumped out of $T$
as a result of this operation. Clearly, one has
$$(T^{\uparrow A})^{\downarrow \eta(T^{\uparrow A}) }=T.
$$

 Direct use of   Garfinkle insertion
algorithm gives the following result where the bold  letters
indicate the domino cell which is pushed back during the reverse
insertion algorithm.

\begin{corollary} \label{reverse.insertion.cor1}
Let $T$ be an $r$-domino tableau and $A$ is a domino corner. Further
let $A'$ be the domino cell which is pushed back by $A$
 in  the first step of the reverse insertion $T^{\uparrow A}$. Then
\begin{enumerate}
\item[$i)$] If  $A=\{(i,j),(i,j+1)\}$  and   $\lb(T,A)=(a,a)$  then
 $A' \subset \{(i-1,k) \mid k\geq j\}$.
\item[$ii)$] If $A=\{(i,j),(i,j+1)\}$  and  $\lb(T,A)=(a',a)$ for some $a'<a$  then
 $A'=\{(i-1,j),(i,j)\}$.

$$\tableau{ {}&{*}& {*}& {*}& {*} \\ {}&{\textbf a}&{\textbf a}\\{}}
                    ~ \longrightarrow~\tableau{{}&{* }& {\textbf *}/a& {\textbf*}/a& {*} \\ {}\\{}}
  ~\mbox{ \ \ \ \ \ \ }~
   \tableau{{}& {*}& {a}& {} \\ {}&{\textbf a'}&{\textbf a}\\{}}
~ \longrightarrow~
\tableau{{}& {\textbf *}/a & {a}& {} \\
                   {}&{\textbf a'}& \\{}}$$

\item[$iii)$] If $A=\{(i,j),(i+1,j)\}$   and   $\lb(T,A)=(a,a)$  then $A'\subset \{(k,j-1)\mid k\geq i\}$.
\item[$iv)$] If $A=\{(i,j),(i+1,j)\}$  and  $\lb(T,A)=(a',a)$ for some $a'<a$ then
$A'=\{(i,j-1),(i,j)\}$.

$$~\tableau{         {}& { }& {} \\
                   {*}&{\textbf{a}}\\{*}&{\textbf{a}}\\{*}}
                    ~ \longrightarrow~
\tableau{         {}& { }& {} \\
                   {*}\\{\textbf *}/a \\{\textbf *}/a }
\mbox{ \ \ \ \ \ \ \  }
~\tableau{         {}& { }& {} \\
                   {*}&{\textbf a'}\\{a}&{\textbf{a}}\\{}}
                    ~ \longrightarrow~
\tableau{         {}& { }& {} \\
                   {\textbf *}/a&{\textbf a'}\\{a}\\{}}
$$
\end{enumerate}
\end{corollary}

\vskip.1in
\begin{example} Let $ T \in SDT^3(5)$  and $B=\{(4,2),(5,3)\}$ as given below. Then one can obtain $T^{\uparrow B}$ in the following manner:

$$   T=\tableau{{ }&{ }&{1}&{1}\\
                   {}&{2}&{2} \\
                   {3}&{3}& {5} \\
                   {4}&{\textbf{4}}&{\textbf{5}}}
                      \rightarrow
\tableau{{ }&{ }&{1}&{1}\\
                   {}&{2}&{2} \\
                   {3}&{\textbf{3}/5}& {5} \\
                   {4}&{\textbf{4}}}
                      \rightarrow
\tableau{{ }&{ }&{1}&{1}\\
                   {}&{2}&{2} \\
                   {\textbf{3}/4}&{\textbf{3}/5}& {5} \\
                   {4}}
                      \rightarrow
\tableau{{ }&{ }&{1}&{1}\\
                   {}&{\textbf{2}/3}&{\textbf{2}/3} \\
                   {4}&{5}& {5} \\
                   {4}}
                      \rightarrow
   \tableau{{ }&{ }&{2}&{2}\\
                   {}&{3}&{3} \\
                   {4}&{5}& {5} \\
                   {4}}  =T^{\uparrow B}              $$
Moreover $\eta(T^{\uparrow B})=1$.
\end{example}

\subsection{ Barbash and Vogan algorithm.}
We will now explain the algorithm which is provided by Barbash and
Vogan in \cite{Barbash-Vogan} to establish the bijection between
signed permutations and standard $r$-domino tableaux for $r=0,1$
where $r=0$ represents type $C$ and $r=1$ represents type $B$ signed
permutations. The extension  of this algorithm for larger cores is
provided in \cite{Geck-Lam}.  We also remark that the equivalence of
Barbash-Vogan algorithm to Garfinkle's algorithm for $r=0,1$, is due
to van Leeuwen~\cite{Leeuwen} .

Recall that for a signed permutation
$\alpha=\alpha_1~\alpha_2\ldots\alpha_n$ the palindrome
representation of $\alpha$ is given by
$\alpha^{0}=\overline{\alpha}_n\ldots
\overline{\alpha}_2~\overline{\alpha}_1~\alpha_1~\alpha_2\ldots\alpha_n$
 where $\overline{\alpha}_i=-\alpha_i$. We
call  $\alpha^{0}$ as {\it $0$-core } representation of $\alpha$.
Clearly $0$-core representation defines an injective map from the
set of all signed permutations of size $n$ into $S_{2n}$.

 By following the approach of
\cite{Geck-Lam} let us describe how to extend this representation
for larger cores. We first  identify $\{1,2,\ldots,r(r+1)/2\}$ with
$\{0_1,0_2,\ldots,0_{r(r+1)/2}\}$ together with the total ordering
$$\bar{n}<\ldots <\bar{2}<\bar{1}<0_1<0_2<\ldots<0_{r(r+1)/2}<1<2\ldots<n.$$
Let $w \in S_{r(r+1)/2}$ be a permutation under this identification,
whose RSK insertion tableau is of shape $(r,r-1,\ldots,1)$. Now for
$\alpha \in B_n$ let {\it $r$-core representation} of $\alpha$ to be
$$ \alpha^{r}=\overline{\alpha}_n\ldots
\overline{\alpha}_2~\overline{\alpha}_1~w~\alpha_1~\alpha_2\ldots\alpha_n.
$$

The algorithm   introduced by Barbash and Vogan for $r=0$ and $r=1$
first applies RSK algorithm  on $\alpha^0$ and respectively
$\alpha^1$.  Then starting from the lowest number $\bar{n}$, it
vacates the negative  integer $\bar{i}$ in the  tableaux by jeu de
taquin slides until it becomes adjacent to $i$, where the evacuation
is repeated for $\overline{i-1}$ until $i=1$. The following example
illustrates this algorithm for $r=1$.

\begin{example} For  $\alpha=3~\bar{1}~2 \in B_n$, we have
 $\alpha^{1}=\bar{2}~1~\bar{3}~0~3~\bar{1}~2$ be its $1$-core representation. Then
Barbash-Vogan algorithm yields:
$$ P(\alpha^{1})=\tableau{\bar{3}&\bar{1}&{2}\\ \bar{2}&0&3\\1}
\mapsto \tableau{\bar{2}&\bar{1}&{2}\\ 0&\bar{3}&3\\1} \mapsto
\tableau{\bar{1}&\bar{2}&2\\ 0&\bar{3}&3\\1}
 \mapsto
\tableau{0&\bar{2}&{2}\\ \bar{1}&\bar{3}&3\\1} \mapsto
\tableau{{}&2&2\\ 1&3&3\\1}=P^1(\alpha).
$$
Similarly $ Q(\alpha^{1})=\tableau{\bar{3}& \bar{2}&1\\
\bar{1}&0&3\\{2}} \mapsto \tableau{{}&1&1\\ 2&3&3\\2}=Q^1(\alpha)$.
\end{example}

\vskip.1in On the other hand by the result of \cite{Geck-Lam}, one
only  needs to apply the same algorithm on $\alpha^r$ in order to
find $r$-domino tableaux $P^r(\alpha)$ and $Q^r(\alpha)$ for larger
cores.

\begin{theorem} [\cite{Geck-Lam}, Theorem 3.3]  \label{Barbash.Vogan.theo1}
Signed permutations $\alpha$ and $\beta$ have the same insertion
$r$-domino tableau if and only if $\alpha^{r}$ and $\beta^{r}$ have
the same RSK insertion tableau.
\end{theorem}

The following two propositions can be deduced by using
Definition~\ref{tableau.def.2}, Theorem~\ref{tableau.theo.0},
Theorem~\ref{tableau.theo.1} and Theorem~\ref{Barbash.Vogan.theo1}.

\begin{proposition}\label{Barbash.Vogan.prop0} Let $\alpha$ be a
signed permutation. Then
$$P^r(\alpha^{-1})=Q^r(\alpha)
~\text{and}~Q^r(\alpha^{-1})=P^r(\alpha).
$$
\end{proposition}

\begin{proposition}\label{Barbash.Vogan.prop1}
Let $\alpha$ and $\beta$ be two signed permutations such that $\alpha^{r}~ \Knuth~ \beta^{r}$.   Then $P^r(\alpha)=P^r(\beta)$, in other
words $\alpha$ and $\beta$ have the same insertion $r$-domino
tableau.
\end{proposition}

\subsection{Descents of domino tableaux and Vogan's map}

Recall that  $B_n$ carries a  Coxeter group structure with the
generator set $S=\{s_0, s_1, \ldots, s_{n-1}\}$ where $\{s_i=(i,i+1)
| 1\leq i\leq n-1\}$ is the set of transpositions which  also
generates the symmetric group $S_n$ and $s_0$ corresponds to the
transposition $(-1,1)$. Let  $l(\alpha)$ denote the length of
$\alpha$, which is the minimum number of generators of $\alpha$ and
let

$$\begin{aligned}
\Des_L(\alpha)&:=\{i ~|~ l(s_i \alpha)< l(\alpha) ~\text{and}~ 0\leq i\leq n-1\} \\
                        &=\{i \mid \text{if }~ 1\leq i \leq n-1~\text{ and }~ i+1
                        ~\text{comes before} ~i ~\text{in}~ \alpha^0\} \cup \{0 \mid \text{ if }~1~\text{comes before}~-1~\text{in}~\alpha^0\}\\
\Des_R(\alpha)&:=\Des_L(\alpha^{-1}) \end{aligned}
$$
 denote respectively the sets of   left and right descents
of $\alpha$.

Now we define the descent set of a $r$-domino tableau $T$ in the
following way:

$$\begin{aligned}\Des(T)&:=\{i \mid \text{ if the domino labeled with
}~ (i+1,i+1) ~\text{ lies below the one labeled with } (i,i) \} \\ &
\cup \{0 \mid \text{ if the domino labeled with }~ (1,1)~ \text{ is
vertical} \}\end{aligned}
$$

It is a well known property of the RSK algorithm that for a
permutation $w\in S_n$, we have
$$\Des_L(w)= \Des(P(w))$$
 where  the descent set
of a  (skew or Young) tableau $T$ is defined by $\Des(T)=\{i \mid
i+1 ~\text{lies below}~ i ~\text{in}~ T \}$.  Now we have:

\begin{proposition}  \label{Barbash.Vogan.prop11} For  $\alpha \in B_n$ we have
$\Des_L(\alpha)=\Des(P^r(\alpha))$.
\end{proposition}
\begin{proof} Observe that $i \in \Des_L(\alpha)$ if and only if one of the following
$$(i+1)i,~(\overline{i+1})i,~i(\overline{i+1}), \bar{i}(\overline{i+1}) $$
is a subsequence in $\alpha$. For the first two cases let $S$
denotes the tableau obtained by inserting all the numbers which
comes before   $i$  in $\alpha$, by  Garfinkle's insertion
algorithm. Therefore  $S$ has a
 domino cell double labeled by $i+1$.  Now since  $i$ is inserted horizontally to the first row of $S$
the domino cell labeled by $i+1$   lies below the one labeled by $i$
in $S^{\downarrow i}$, and moreover it remains to be below until the
last letter in $\alpha$ inserted, since $i$ and $i+1$ are
consecutive numbers. For the last two cases let $T$ be the tableau
obtained by inserting all numbers which comes before $i+1$ in
$\alpha$. This time domino  cell double labeled by $i$ lies in $T$
and since $i+1$ is inserted vertically to the first column of the
tableau of $T$, this vertical cell lies below the one labeled by $i$
in $T^{\downarrow i}$. On the other hand insertion of subsequent
numbers in $\alpha$ does not change this rule and hence $i\in
\mathrm{Des}(P^r(\alpha))$.

For the reverse inclusion observe that $i \not\in \Des_L(\alpha)$ if and only if one of the following
$$i(i+1), ~ \bar{i}({i+1}),~(i+1)\bar{i},~(\overline{i+1})\bar{i} $$
is a subsequence in $\alpha$. In this case a similar argument to the
one used above shows that  $i \not\in \Des(P^r(\alpha))$.

  \end{proof}

\begin{figure}[h]
\vspace{-1.8in}
\begin{center}
$\begin{array}{c} \epsfysize=7in \epsffile{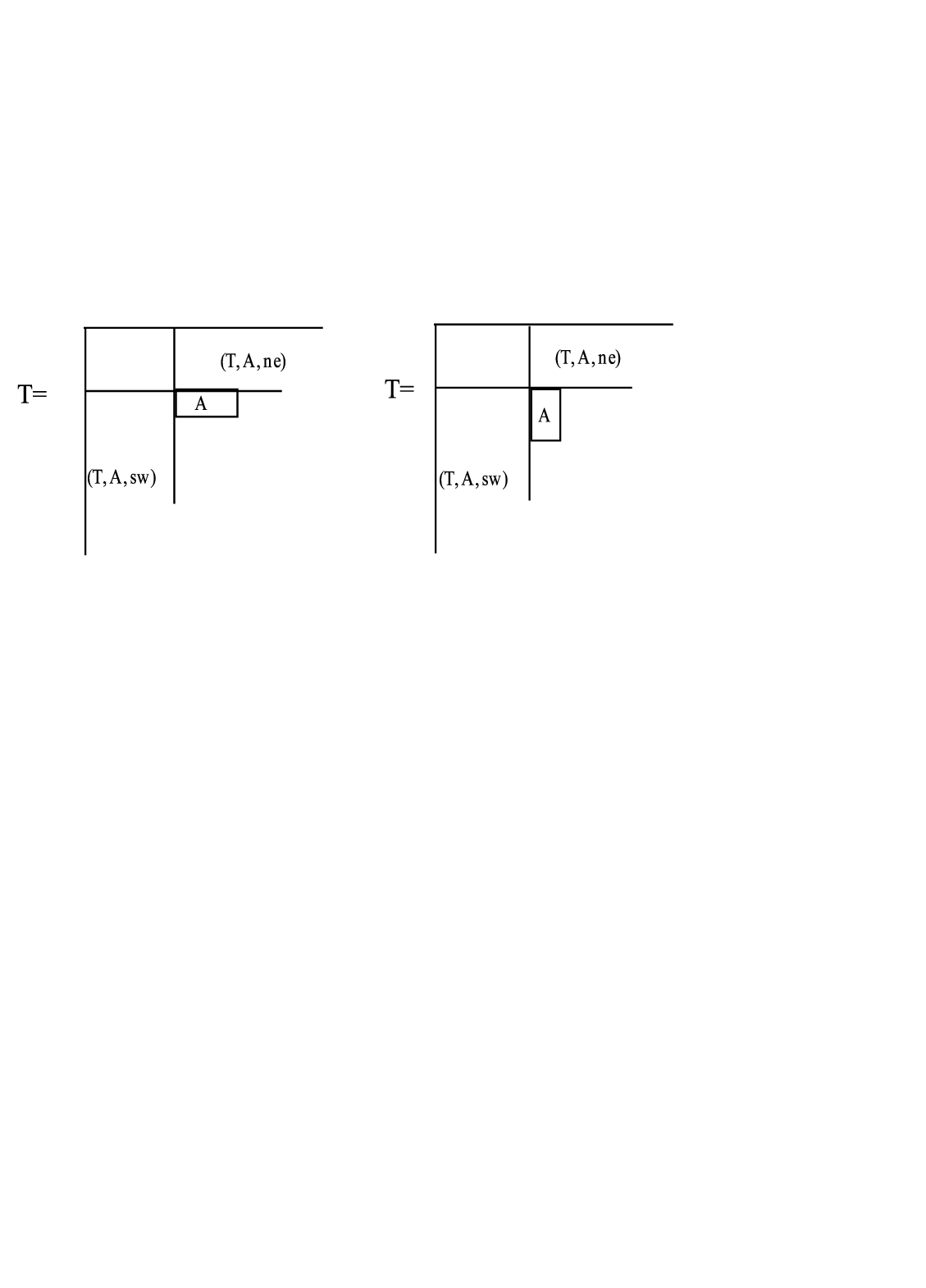}
\end{array}$
\end{center} \vspace{-4in}
\caption{ } \label{reverse.insertion.fig1}
\end{figure}

\begin{definition} Let $T$ be a $r$-domino tableau and $A$ be a  domino corner of  $\sh(T)$  such that
$A=\{(i,j),(i,j+1)\}$ or $A=\{(i,j),(i+1,j)\}$. We denote by
$(T,A,\mathrm{ne})$ and $(T,A,\mathrm{ne})$  the regions of $T$ such
that
$$\begin{aligned}(T,A,\mathrm{ne}):=&\{ (k,l) ~\mid~  k<i ~\text{and}~ l\geq j
\}\\
(T,A,\mathrm{sw}):=&\{ (k,l) ~\mid~  k\geq i ~\text{and}~ l<j \}
\end{aligned}
$$
as illustrated in Figure~\ref{reverse.insertion.fig1}.
\end{definition}

\input{epsf}

\vskip.1in
 Now we are ready to give the following lemma which is crucial in the proof of
 Theorem~\ref{main.theorem2}.

\begin{lemma} \label{reverse.insertion.lem2}
 Let $T$ be a $r$-domino tableau and $A$ be  a domino corner of $\sh(T)$.
\begin{enumerate}
\item[$i)$] Suppose  $B$  is a domino corner  of $\sh(T^{\uparrow A})$ which
lies in the portion  $(T, A,\mathrm{sw})$. Then
$$ \eta( T^{\uparrow A\uparrow  B} )< \eta(T^{\uparrow A}). $$
\item[$ii)$]  Suppose $B$ is a  domino corner of $\sh(T^{\uparrow A})$ which
 lies in the portion  $(T, A,\mathrm{ne})$. Then
$$ \eta( T^{\uparrow A \uparrow B}) >  \eta(T^{\uparrow A}). $$
\end{enumerate}
\end{lemma}

\begin{proof} We will just prove the first part of the theorem since the same method applies to the second part.
Let $a=\eta(T^{\uparrow A})$, $b=\eta( T^{\uparrow A\uparrow  B})$
and $u$ be a word such that $P^r(u)=T^{\uparrow A\uparrow  B}$. Then
clearly the sign permutation $\alpha=uba$ has
$$P^r(\alpha)=P^r(uba)=P^r(u)^{\downarrow b\downarrow a}=(T^{\uparrow A\uparrow  B})^{\downarrow b\downarrow
a}=(T^{\uparrow A})^{\downarrow a}=T
$$
and
$$\lb(Q^r(\alpha),A)=(n,n) ~\text{and}~\lb(Q^r(\alpha),B)=(n-1,n-1).$$
On the other hand since $B \in (Q^r(\alpha), A,\mathrm{sw})$ this
shows that $n-1\not \in \Des(Q^r(\alpha))$. Now by
Proposition~\ref{Barbash.Vogan.prop0} and
Proposition~\ref{Barbash.Vogan.prop11} we have
$$n-1\not \in
\Des_L(\alpha^{-1})=\Des_R(\alpha)$$ and therefore $\alpha_{n-1}=b$
can not be bigger than $\alpha_{n}=a$. Therefore $a=\eta(T^{\uparrow
A})> b=\eta( T^{\uparrow A\uparrow  B})$ as desired.
\end{proof}

\subsubsection{Vogan's map} Let  $\alpha,\beta \in B_n$ whose usual representations satisfy
  $$\begin{aligned}\alpha=&\alpha_1 \ldots \alpha_{i-1} (\alpha_{i} ~\alpha_{i+1}) ~\alpha_{i+2} \ldots
  \alpha_n\\
 \beta=& \alpha_1 \ldots \alpha_{i-1} (\alpha_{i+1} ~\alpha_{i}) ~\alpha_{i+2} \ldots
 \alpha_n \end{aligned}$$
where  either $ \alpha_i <\alpha_{i+2}<\alpha_{i+1}$ or  $
\alpha_{i} <\alpha_{i-1}<  \alpha_{i+1}$. By extending
 Definition~\ref{tableau.def.2},  we say $\alpha$ and
$\beta$ are  equivalent through single Knuth relation.   We write
$\alpha~\Knuth~\beta$  in $B_n$ if one of them can be obtained from
the other by applying a sequence of Knuth relations. Clearly
$\alpha~\Knuth~\beta$ yields  $\alpha^r~ \Knuth~ \beta^r$ i.e.,
Knuth relations which are obtained on the usual representation of
signed permutations does not change insertion $r$-domino tableaux.
On the other hand their effect on recording tableaux for the case
$r=0,1$ is studied by Garfinkle \cite[2.1.10--2.1.19]{Garfinkle2}.
Our following analysis is based on her work, including all
notations, definitions and  maps.  Let $r\geq 0$.

For $i,j$ two adjacent integers satisfying  $1\leq i, j \leq n-1$,
consider the following sets:
$$\begin{aligned} D_{i,j}(B_n):=&\{\alpha \in B_n \mid i  \in \Des_L(\alpha)~\text{ but }~j \not \in \Des_L(\alpha)\}\\
D_{i,j}(SDT^r(n)):=&\{T\in SDT^r(n) \mid i \in \Des(T)~\text{ but
}~j \not \in \Des(T)\}
\end{aligned}$$
together with the map  $V_{i,j}: D_{i,j}(B_n) \mapsto D_{j,i}(B_n)$
where $V_{i,j}(\alpha)=\{s_i\cdot\alpha, ~s_{j}\cdot\alpha\} \cap
D_{j,i}(B_n).$ Also  define a map $$V_{i,j}: D_{i,j}(SDT^r(n))
\mapsto D_{j,i}(SDT^r(n))$$
 in the following manner: Without loss of generality we assume that $j>i$, i.e.,  $j=i+1$.
Observe that if $i \in\Des(T)$ but $i+1 \not\in \Des(T)$ then $i+1$
lies  below $i$ in $T$ whereas $i+2$ lies right to
$i+1$ in $T$. On the other hand we have two cases according to the
positions of dominos labeled with $(i,i)$ and $(i+2,i+2)$ with
respect to each other.

 \noindent {\bf Case 1}.  We first assume that   $i+2$ lies below $i$
in $T$. Since the   $i+2$ lies  to the right of
 $i+1$ and $i+1$ lies below $i$ we have two cases to consider: If the boundaries  $\mathrm{Dom}(T,i+1)$
 and $\mathrm{Dom}(T,i)$  intersect at most at a point then $V_{i,i+1}(T)$ is
obtained by interchanging the labels $i$ and $i+1$ in $T$. Otherwise
there is only one possibility  which satisfies  $i+2$ lies below $i$
and it lies to the right of $i+1$,  in which  $T$ has the subtableau
$U$ as illustrated below and $V_{i,i+1}(T)$ is obtained by
substituting $U$ with $U'$ in $T$.
 $$U= \tableau{_{i}&_{i}\\_{i+1}&_{i+2}\\_{i+1}&_{i+2}}~\mbox{ \ }~
U'= \tableau{_{i}&_{i+1}\\_{i}&_{i+1}\\_{i+2}&_{i+2}}
$$

\noindent {\bf Case 2}. Now we assume  $i+2$ lies strictly right to
$i$ in $T$. Again
 if the boundaries of $\mathrm{Dom}(T,i+1)$ and $\mathrm{Dom}(T,i+2)$
intersect at most at a point then $V_{i,i+1}(T)$ is obtained by
interchanging the labels $i+1$ and $i+2$ in $T$. Otherwise there is
only  one possible case where  $T$ has the subtableau $U$ given
below  and $V_{i,i+1}(T)$ is obtained by substituting $U$ with $U'$
in $T$.
$$U=\tableau{_{i}&_{i}&_{i+2}\\_{i+1}&_{i+1}&_{i+2}}~\mbox{ \ }~
U'=\tableau{_{i}&_{i+1}&_{i+1}\\_{i}&_{i+2}&_{i+2}}
$$

\begin{example} We have $T_2=V_{5,6}(T_1)$, $T_3=V_{3,4}(T_2)$, and
$T_4=V_{4,5}(T_3)=V_{6,5}(T_3)$  for the following tableaux.
$$T_1=\tableau{{1}&{2}&{\mathbf 5}\\{1}&{2}&{\mathbf 5}\\{3}&{3}&{\mathbf 7}\\{4}&{\mathbf 6}&{\mathbf 7}\\{4}&{\mathbf 6}}
~\mbox{ }~ T_2=\tableau{{1}&{2}&{6}\\{1}&{2}&{6}\\{\mathbf
3}&{\mathbf 3}&{7}\\{\mathbf 4}&{\mathbf 5}&{7}\\{\mathbf
4}&{\mathbf 5}} ~\mbox{ }~ T_3=\tableau{{1}&{2}&{\mathbf
6}\\{1}&{2}&{\mathbf 6}\\{3}&{\mathbf 4}&{\mathbf 7}\\{3}&{\mathbf
4}&{\mathbf 7}\\{\mathbf 5}&{\mathbf 5}} ~\mbox{ }~
T_4=\tableau{{1}&{2}&{5}\\{1}&{2}&{5}\\{3}&{4}&{7}\\{3}&{4}&{7}\\{6}&{6}}
$$
\end{example}
\begin{remark} The map $V_{i,j}$ is first introduced  on the symmetric group by  Vogan  \cite{Vogan},
 with the aim of classifying the primitive ideals in the universal enveloping algebra of complex semi simple Lie algebras.
\end{remark}

\begin{lemma} \label{Barbash.Vogan.lem1}
Let $i$ and $j$ be two consecutive  integers such that  $1\leq i,j
\leq n-1$.  Suppose $\alpha \in D_{i,j}(B_n)$. Then $ P^r(\alpha)
\in D_{i,j}(SDT^r(n))$ and
$$P^r( V_{i,j}(\alpha))= V_{i,j}(P^r(\alpha)).
$$
\end{lemma}
\begin{proof} This result is first proven by Garfinkle  \cite[Theorem~2.1.19]{Garfinkle2} for  $r=0,1$.
On the other hand one can check that her proof does not depend  on
the specific  value of $r$ and  it can easily be extended for any
value of $r$. We omit the proof for the sake of space.
\end{proof}

The following result has an important role in the proof
of Theorem~\ref{main.theorem2}.

\begin{corollary} \label{Barbash.Vogan.cor2}
Suppose $\alpha=\alpha_1\ldots \alpha_{i-1}(\alpha_i
\alpha_{i+1})\alpha_{i+2}\ldots \alpha_n $ and $
\beta=\alpha_1\ldots \alpha_{i-1}(\alpha_{i+1}
\alpha_{i})\alpha_{i+2}\ldots \alpha_n $ differ by a single
Knuth relation. Then one of the following is satisfied:
\begin{enumerate}
\item[1)] $\alpha_i<\alpha_{i+2}<\alpha_{i+1}$ then $\beta^{-1}=V_{i+1,i}(\alpha^{-1})$ and
$Q^r(\beta)=V_{i+1,i}(Q^r(\alpha))$.
\item[2)] $\alpha_i>\alpha_{i+2}>\alpha_{i+1}$ then  $\beta^{-1}=V_{i,i+1}(\alpha^{-1})$ and
$Q^r(\beta)=V_{i,i+1}(Q^r(\alpha))$.
\item[3)] $\alpha_i<\alpha_{i-1}<\alpha_{i+1}$ then $\beta^{-1}=V_{i-1,i}(\alpha^{-1})$ and $
Q^r(\beta)=V_{i-1,i}(Q^r(\alpha))$.
\item[4)] $\alpha_i>\alpha_{i-1}>\alpha_{i+1}$ then
 $\beta^{-1}=V_{i,i-1}(\alpha^{-1})$ and $
Q^r(\beta)=V_{i,i-1}(Q^r(\alpha))$.
\end{enumerate}

\end{corollary}
\begin{proof} Assumptions on $\alpha$ and $\beta$ in the first case  yields that
  $\alpha^{-1} \in D_{i+1,i}(B_n)$ and  $\beta^{-1}=s_i \cdot \alpha^{-1} \in D_{i,i+1}(B_n)$. Therefore
  $\beta^{-1}=V_{i+1,i}(\alpha^{-1})$ and by Lemma~\ref{Barbash.Vogan.lem1}
$$Q^r(\beta)=P^r(\beta^{-1})=P^{r}(V_{i+1,i}(\alpha^{-1}))=V_{i+1,i}(P^r(\alpha^{-1}))=V_{i+1,i}(Q^r(\alpha)).$$
  For the other cases the result follows similarly.
\end{proof}


\section{ Plactic relations for $r$-domino tableaux}
\label{proof.of.main.theo.sec} \noindent

Recall that for any $a\in\mathbb{Z}$, $\bar{a}$ represents $-a$ if
$a>0$ and it represents $|a|$ otherwise.

\begin{definition}\label{plactic.relations}

For  $\alpha$ and $ \beta$ are two signed permutations  in $B_m$ and
$r\geq0$, we say  $\alpha$ and $\beta$ are {\it $r$-plactic
equivalent}, $\alpha~\rp~\beta$, if one of them can  be obtained
from the other by applying a sequence of $\mathrm{D}_i^r$ relations
for $i=1,\ldots 5$,  explained below. Moreover, we say $\alpha$ and
$\beta$ are {\it $r$-coplactic equivalent}, $\alpha~\rcp~\beta$, if
$\alpha^{-1}~\rp~\beta^{-1}$.

Let $\alpha=\alpha_1 \ldots  \alpha_{m} \in B_m$.

\begin{enumerate}
\item[$\mathrm{D}_1^r$:] If  $ \alpha_i <\alpha_{i+2}<  \alpha_{i+1}$ or $ \alpha_i <\alpha_{i-1}<  \alpha_{i+1}$
for some $i\leq m-1$, then
$$\alpha=
\alpha_1 \ldots \alpha_{i-1}~(\alpha_{i} ~\alpha_{i+1})
~\alpha_{i+2} \ldots \alpha_m~  \sim~
 \alpha_1 \ldots \alpha_{i-1}~ (\alpha_{i+1} ~\alpha_{i}) ~\alpha_{i+2} \ldots
 \alpha_m$$
\item[$\mathrm{D}_2^r$:] If $r\geq1$ and if there exists  $0<j\leq r$ such that $\alpha_j$ and $\alpha_{j+1}$ have opposite signs   then
$$\alpha=\alpha_1 \ldots  \alpha_{j-1}(\alpha_{j} ~\alpha_{j+1}) \ldots \alpha_{r+2}\ldots  \alpha_m ~\sim~
\alpha_1 \ldots  (\alpha_{j+1} ~\alpha_{j}) \ldots \alpha_{r+2}
\ldots\alpha_m $$
\item[$\mathrm{D}_3^r$:] Suppose that  $|\alpha_1|> |\alpha_i|$ for all $2 \leq i \leq r+2$
and  $\alpha_2 \ldots \alpha_{r+2}$ is obtained by  concatenating
some positive decreasing sequence to the end of some negative
increasing sequence (or vice versa), where at least one of the
sequences is nonempty. Then
$$\alpha=\alpha_1~\alpha_2 \ldots \alpha_{r+2} \ldots
\alpha_m ~\sim ~ \overline{\alpha_1}~\alpha_2 \ldots \alpha_{r+2}\ldots  \alpha_m $$
\item[$\mathrm{D}_4^r$:]  Let
for some $k\geq 1$, $s=(k+1)(r+k+1)\leq m $ and $u= \alpha_1\ldots
\alpha_{s-1}$ is obtained by concatenating the sequences
$a_{i,i+r}\ldots a_{i,1}$ and $b_{i,i}\ldots b_{i,1}$ for $1\leq i
\leq k$  and $a_{k+1,r+k}\ldots a_{1,1}$  in the following manner:
$$
u= \alpha_1\ldots \alpha_{s-1}=a_{_{1,r+1}}\ldots
a_{_{1,1}}b_{_{1,1}}\ldots a_{_{k,r+k}}\ldots a_{_{k,1}}
b_{_{k,k}}\ldots b_{_{k,1}}a_{_{k+1,r+k}}\ldots a_{_{k+1,1}}
$$
where  the integers $a_{_{i,j}}$ and $b_{_{i,j}}$, if exist in $u$,
satisfy the following conditions:
$$\begin{aligned}
&a_{_{i,j}}>0 ~\mbox { and }~b_{_{i,j}}<0 ~\mbox{(or vice versa)}\\
&|a_{_{i,j-1}}|<|a_{_{i,j}}|<|a_{_{i+1,j}}| ~\mbox{ and }~  |b_{_{i,j-1}}|<|b_{_{i,j}}|<|b_{_{i+1,j}}|\\
& |b_{_{i,i}}|<|a_{_{i+1,r+i+1}}|<|b_{_{i+1,i+1}}|~\mbox{for all
}~i=1,\ldots, k-1.
\end{aligned}$$

Let $n=\mathrm{max}\{ |\alpha_1|,\ldots,|\alpha_{s-1}| \}$ and
suppose that the integer $\alpha_{s}=z$ satisfies one of the
followings:

\begin{enumerate}
\item[$i.$] $|b_{_{k,k}}|=n$ and  $z$ is an integer between $a_{_{k+1,1}}$ and
$b_{_{k,1}}$
\item[$ii.$] $|a_{_{k+1,r+k}}|=n$ and  $z$ is  an integer between
$a_{_{k,1}}$ and $b_{_{k,1}}$
\item[$iii.$] $|a_{_{k+1,r+k}}|=n$, $z$ is  an integer between
$a_{_{k,1}}$ and $a_{_{k+1,1}}$ and $|a_{_{k+1,i}}|<|a_{_{k,i+1}}|$
\end{enumerate}
 Then we set $$\alpha=uz\alpha_{s+1}\ldots \alpha_m
~\sim ~u'z\alpha_{s+1}\ldots \alpha_m$$ where
$u'=a_{_{1,r+1}}\ldots a_{_{1,1}}b_{_{1,1}} \ldots
\overline{b_{_{k,k}}}a_{_{k,r+k}}\ldots a_{_{k,1}}b_{_{k,k-1}}\ldots
b_{_{k,1}}a_{_{k+1,r+k}}\ldots a_{_{k+1,1}}$.\\

\item[$\mathrm{D}_5^r$:]  Let
for some $k\geq 1$, $s=(k+1)(r+k+2)\leq m $ and $u= \alpha_1\ldots \alpha_{s-1}$ is
 obtained by concatenating the sequences    $a_{i,i+r}\ldots a_{i,1}$ and $b_{i,i}\ldots b_{i,1}$ for $1\leq i \leq k$   and $a_{k+1,r+k+1}\ldots a_{1,1}$  and $b_{k+k,k}\ldots a_{1,1}$ in the following manner:
$$
u= \alpha_1\ldots \alpha_{s-1}=a_{_{1,r+1}}\ldots
a_{_{1,1}}b_{_{1,1}}\ldots a_{_{k,r+k}}\ldots a_{_{k,1}}
b_{_{k,k}}\ldots b_{_{k,1}}a_{_{k+1,r+k+1}}\ldots
a_{_{k+1,1}}b_{_{k+1,k}}\ldots b_{_{k+1,1}}
$$
where  the integers $a_{_{i,j}}$ and $b_{_{i,j}}$, if exist in $u$,
satisfy the followings:
$$\begin{aligned}
&a_{_{i,j}}>0 ~\mbox { and }~b_{_{i,j}}<0 ~\mbox{(or vice versa)}\\
&|a_{_{i,j-1}}|<|a_{_{i,j}}|<|a_{_{i+1,j}}| ~\mbox{ and }~  |b_{_{i,j-1}}|<|b_{_{i,j}}|<|b_{_{i+1,j}}|\\
& |a_{_{i,r+i}}|<|b_{_{i,i}}|<|a_{_{i+1,r+i+1}}|~\mbox{for all
}~i=1,\ldots, k.
\end{aligned}$$

Let $n=\mathrm{max}\{ |\alpha_1|,\ldots,|\alpha_{s-1}| \}$ and
suppose that the integer $\alpha_{s}=z$ satisfies one of the
followings:
\begin{enumerate}
\item[$i.$] $|a_{_{k+1,r+k+1}}|=n$ and  $z$ is an integer between
$a_{_{k+1,1}}$ and $b_{_{k+1,1}}$
\item[$ii.$] $|b_{_{k+1,k}}|=n$ and  $z$ is an integer between
$a_{_{k+1,1}}$ and $b_{_{k,1}}$
\item[$iii.$] $|b_{_{k+1,k}}|=n$, $z$ is an integer between
$b_{_{k,1}}$ and $b_{_{k+1,1}}$ and  $|b_{_{k+1,i}}|<|b_{_{k,i+1}}|$
for some $1<i\leq k-1$.
\end{enumerate}
 Then we set $$\alpha=uz\alpha_{s+1}\ldots \alpha_m
~\sim ~u'z\alpha_{s+1}\ldots \alpha_m$$ where
$u'=_{_{1,r+1}}\ldots a_{_{1,1}}b_{_{1,1}}\ldots a_{_{k,r+k}}\ldots a_{_{k,1}}
\overline{a_{_{k+1,r+k+1}}}b_{_{k,k}}\ldots
b_{_{k,1}}a_{_{k+1,r+k}}\ldots a_{_{k+1,1}} b_{_{k+1,k}}\ldots
b_{_{k+1,1}}$.

\end{enumerate}
\end{definition}

Following example illustrates $D_4^r$ for $r=0$, where $k=1$.
\begin{example}  Consider
$$T=\tableau{{1}&{1}&{3}&{3}\\{2}&{4}\\{2}&{4}} \in SDT^0(4)$$ and its
 domino corners
$A=\{(3,1),(3,2)\},B=\{(2,2),(3,2)\},C=\{(1,3),(1,4)\} $ together
with the sets
$$ \mathcal{K}_A(T)=\{143\bar{2}, 413\bar{2},  \bar{4}13\bar{2}\},~
 \mathcal{K}_B(T)=\{ 13 \bar{4}\bar{2}, 1\bar{4}3\bar{2}\},~\mathcal{K}_C(T)=\{\bar{4}1\bar{2}3,
 1\bar{4}\bar{2}3 ,41\bar{2}3 \}
$$
where $\mathcal{K}_{(-)}(T)$  consist of all signed permutations
whose insertion gives $T$ with the rule that the last opening domino
corner cell in  the insertion is $(-)$. One can easily see that
permutations in each set above   are related by a sequence of
$D_1^0$ and  $D_3^0$ relations. Therefore one needs to obtain some
relations between the permutations of  these distinct sets. Here $
413\bar{2} \in \mathcal{K}_A(T)$ and  $41\bar{2}3\in
\mathcal{K}_C(T)$ are related by $D_1^r$ relations so  a relation
which connects either $\mathcal{K}_A(T)$ and  $\mathcal{K}_B(T)$ or
$\mathcal{K}_B(T) $ and $\mathcal{K}_C(T)$ is needed. Now  $D_4^0$
relates $1\bar{4}3\bar{2}\in\mathcal{K}_B $ and
$413\bar{2}\in\mathcal{K}_A(T)$  under the formulation
$$\begin{aligned} &1\bar{4}3\bar{2}=a_{11}b_{11}a_{21}z\\
&413\bar{2}=\overline{b_{11}}a_{11}a_{21}z.\end{aligned}$$ where
$k=1$ and $r=0$.
\end{example}

The following remark will be much more clear in the  proof (Case
3.2.2) of our main result Theorem~\ref{main.theorem2}.

\begin{remark}In general consider the tableau  $T$ whose
shape has the form
$$(s+i,s+i-1,\ldots,s+1,\underline{s,s-2,s-2,s-4},s-3,\ldots,2,1) ~\mbox{for some }~ s\geq 4,
i\geq 0
$$
as Figure~\ref{main.theorem.fig13}(a) and (b) illustrates.
 Observe that such a tableau has
exactly three domino corners   namely
$$A=\{(i+3,s-3),(i+3,s-2)\},B=\{(i+2,s-2),
(i+3,s-2)\},C=\{(i+1,s-1),(i+1,s)\}
$$
Now if $s+i=r+2k$ for some $k\geq 1$ then $\mathcal{K}_A(T)$ and
$\mathcal{K}_B(T)$ are related by  $D_4^r$ relations with
$a_{_{i,j}}>0$ and  $b_{_{i,j}}<0$. On the other hand if
$s+i=r+2k+1$ for some $k\geq 1$ then $\mathcal{K}_A(T)$ and
$\mathcal{K}_B(T)$ are related by  $D_5^r$ relations with
$a_{_{i,j}}<0$ and  $b_{_{i,j}}>0$.

\input{epsf}
\begin{figure}[h]
\vspace{-.6in}
\begin{center}
$\begin{array}{c} \epsfysize=7in \epsffile{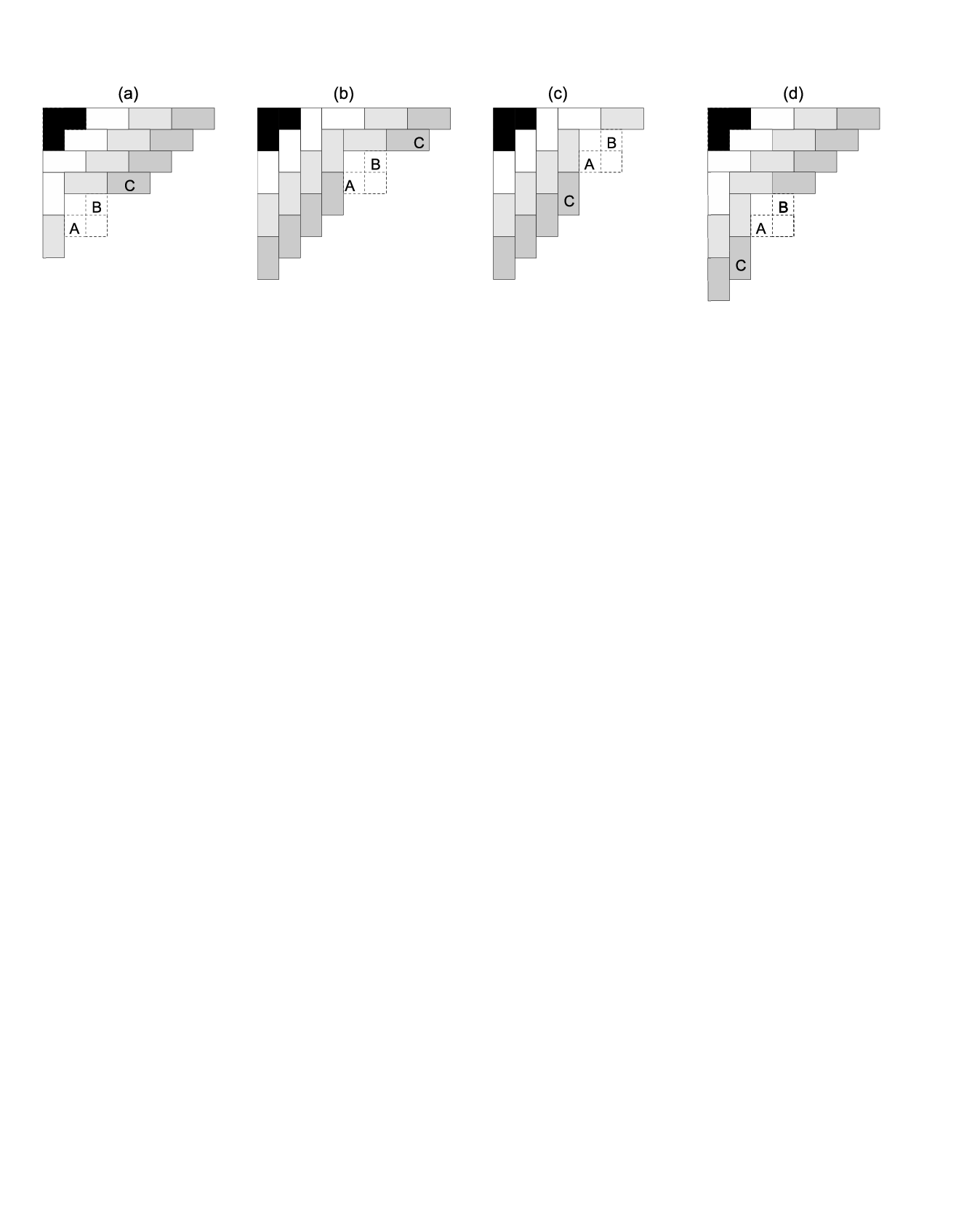}
\end{array}$
\end{center} \vspace{-5.3in}
\caption{ }\label{main.theorem.fig13}
\end{figure}

On the other hand if  $T$ has   shape
$$(s+i,s+i-1\ldots,
s+1,\underline{s,s,s-2,s-2},s-3,\ldots,2,1) ~\mbox{for some }~ s\geq
3, i\geq 0,$$ as Figure~\ref{main.theorem.fig13}(c) and (d)
illustrates, then $T$  has  exactly three domino corners
$$A=\{(i+2,s-1),(i+2,s)\},B=\{(i+1,s),
(i+2,s)\},C=\{(i+3,s-2),(i+3,s-2)\}
$$
Now if $s+i=r+2k$ for some $k\geq 1$ then $\mathcal{K}_A(T)$ and
$\mathcal{K}_B(T)$ are related by  $D_5^r$ relations with
$a_{_{i,j}}>0$ and  $b_{_{i,j}}<0$. On the other hand if
$s+i=r+2k+1$ for some $k\geq 1$ then $\mathcal{K}_A(T)$ and
$\mathcal{K}_B(T)$ are related by  $D_4^r$ relations with
$a_{_{i,j}}<0$ and $b_{_{i,j}}>0$.

Here it is natural to ask whether  some simpler relations exist for
the tableaux presented above. In fact we know  a rule which relates
$\mathcal{K}_B(T)$ and $\mathcal{K}_C(T)$, but its formal
description requires four relations, yet they appear  to us not
handy when it comes to proving our main result.
\end{remark}

 \begin{theorem}\label{main.theorem1} If $\alpha~\rp~\beta$ in $B_m$
  then  they have the same insertion $r$-domino tableaux.
 \end{theorem}

\begin{proof} For the proof of the theorem it will  be enough to consider the case when
 $\alpha$ and $\beta$ differ by a single $\mathrm{D}^r_i$ relation for
 $i=1,\dots,5$.

 Observe that  $\mathrm{D}^r_1$ is just the single Knuth relation defined
 on the usual representation of signed permutations
 and therefore  result follows from
 Proposition~\ref{Barbash.Vogan.prop1}.

  For
 $\mathrm{D}^r_2$, $r\geq1$, let  $\alpha_j$ and $\alpha_{j+1}$ have opposite signs
 in
$\alpha=\alpha_1 \ldots  \alpha_{j-1}(\alpha_{j} ~\alpha_{j+1})
\ldots \alpha_{r+2}\ldots \alpha_m$ for some $j\leq r$. Observe that
the size of empty corners of $r$-staircase shape  is $r+1$.  Let $S$
denote the tableau obtained by inserting  first $j-1$ elements of
$\alpha$ in to $r$-staircase shape. If $r=1$ then $\alpha_1 \ldots
\alpha_{j-1}$ is empty and one can check easily that insertion of
$\alpha_j\alpha_{j+1}$ and $\alpha_{j+1}\alpha_{j}$ into
$1$-staircase shape creates the same tableau.  For $r>1$, observe
that  since $j+1\leq r+1$, the insertion of $\alpha_1 \ldots
\alpha_{j}\alpha_{j+1}$ creates two connected union  of domino
cells, where one consists of horizontal domino cells labeled by the
positive numbers in $\alpha_1 \ldots
 \alpha_{j+1}$,  concatenated to the right of $r$-staircase shape and the other consists of vertical domino
 cells labeled by the absolute value of negative numbers in the same
 sequence, concatenated below  $r$-staircase shape. W.L.O.G we
 assume $\alpha_{j}>0$ and $\alpha_{j+1}<0$ and let $A$ and $B$
 denote the horizontal and vertical  domino cells appearing after the insertion of $\alpha_{j}$ and
 $\alpha_{j+1}$ respectively as illustrated in
 Figure~\ref{proof51} for $r=6$ and  $j=6$.
\input{epsf}
\begin{figure}[h]
\vspace{-.4in}
\begin{center}
$\begin{array}{c} \epsfysize=6in \epsffile{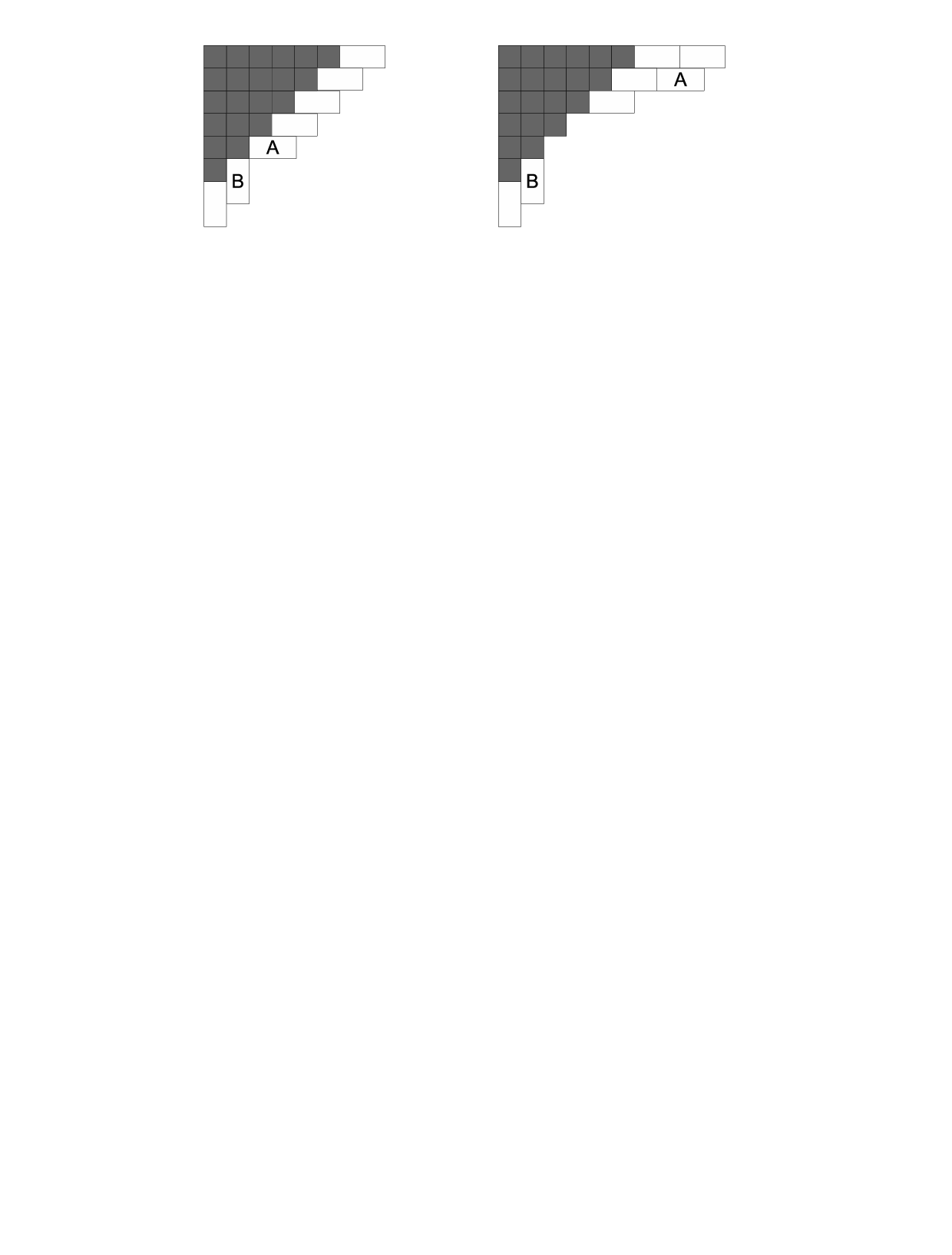}
\end{array}$
\end{center} \vspace{-5in}
\caption{ }\label{proof51}
\end{figure}
Now one can easily  see that changing the order of
  $\alpha_j\alpha_{j+1}$ in $\alpha$ still gives  the same
  tableau,  since $\alpha_j$ and $\alpha_{j+1}$ have opposite
  signs and $j+1\leq r+1$.

 For $\mathrm{D}^r_3$,  let $|\alpha_1|>|\alpha_i|$ for all $2\leq i\leq
 r+2$
 and suppose that for some $k\geq 0$, $l\geq 0$ and $k+l=r+1$
 $$\alpha_2 \ldots \alpha_{r+2}=x_1\ldots x_k~y_1\ldots y_l$$ where
 $x_1\ldots x_k$ is a positive decreasing and  $y_1\ldots y_l$ is a negative increasing
  sequence  (or vice versa). Here observe that we can not have both $k=0$ and $l=0$, since
  then $x_1\ldots x_k~y_1\ldots y_l$
is
  empty but
$\alpha_2 \ldots \alpha_{r+2}$ is  not,  even if $r=0$. So W.L.O.G.
we assume that $l\geq 1$. If $r=0$, we must have $k=0$,  $l=1$ and
$\alpha_2=y_1$ and one can easily check that in this case
$\alpha_1\alpha_2$ and $\overline{\alpha_1}\alpha_2$ give the same
tableau. For $r>0$ the insertion of $\alpha_1 x_1\ldots
x_k~y_1\ldots y_{l-1}$ and $\overline{\alpha_1} x_1\ldots
x_k~y_1\ldots y_{l-1}$ yields two tableaux which differ by only the
position of the domino cell $\{(k+1,l),(k+1,l+1)\}$ and
$\{(k+1,l),(k+2,l)\}$, labeled by $(|\alpha_1|,|\alpha_1|)$, as
illustrated in Figure~\ref{proof41} for $r=5$. On the other hand the
insertion of $y_l$ in both tableaux yields the same tableau.
\input{epsf}
\begin{figure}[h]
\vspace{-.4in}
\begin{center}
$\begin{array}{c} \epsfysize=6in \epsffile{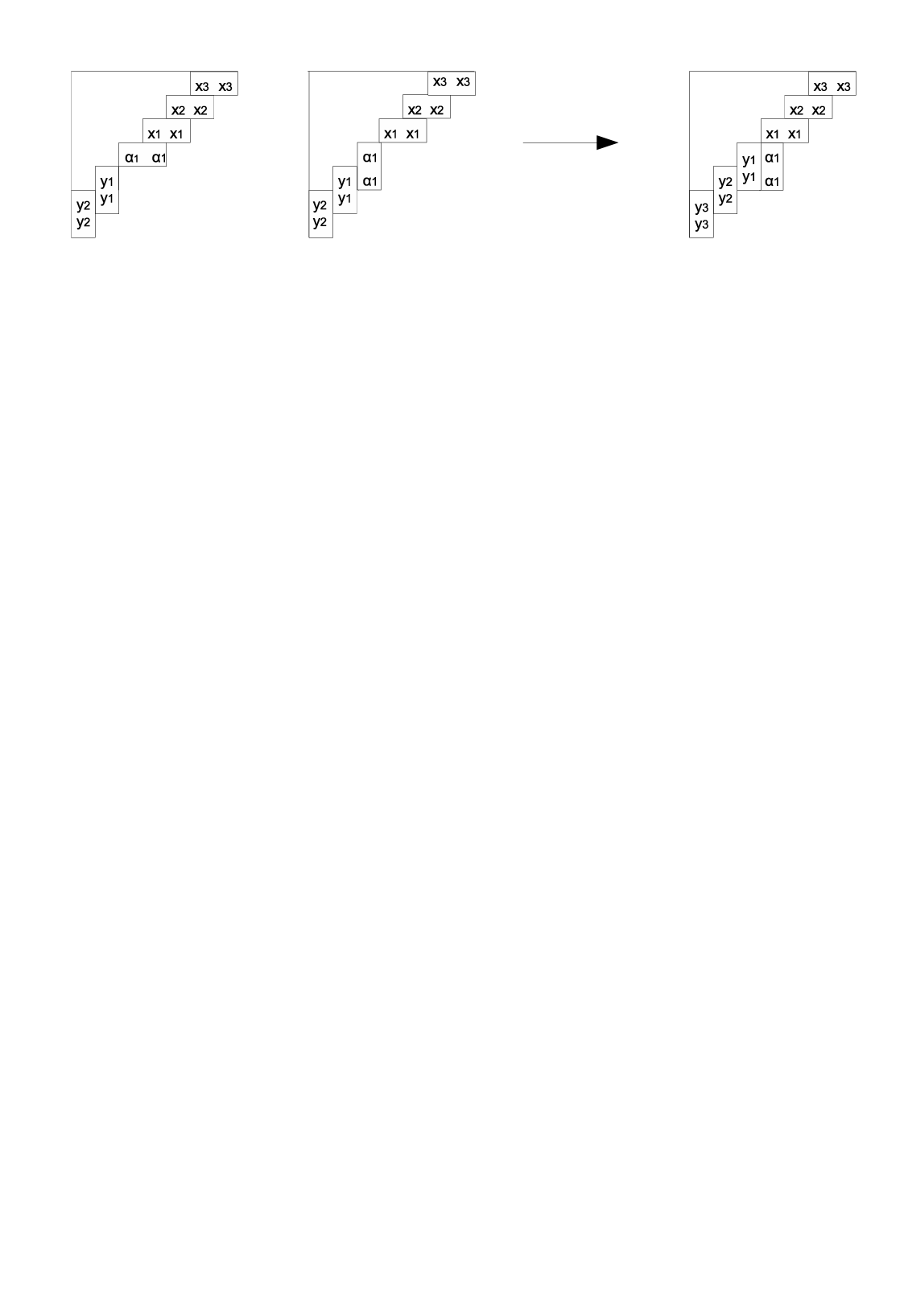}
\end{array}$
\end{center} \vspace{-5.2in}
\caption{ }\label{proof41}
\end{figure}

 In the following we will just deal with the relation $\mathrm{D}^r_4$
 since then  the same method also applies to the relation $\mathrm{D}^r_5$.

It is enough to consider the case  where $\alpha=uz$ and $\beta=u'z$
where $u,u'$ and $z$ are as described in  $\mathrm{D}^r_4$. Observe
that  in this case we have  $m=s=n$.  On the other hand, since the
tableau obtained by taking all $a_{i,j}$'s negative and $b_{i,j}$
positive is the transpose of the tableau obtained by otherwise, it
is enough to   consider only one case. So let
$$\begin{aligned}
u&= a_{_{1,r+1}}\ldots a_{_{1,1}}b_{_{1,1}}\ldots a_{_{k,r+k}}\ldots
a_{_{k,1}} b_{_{k,k}}\ldots b_{_{k,1}}a_{_{k+1,r+k}}\ldots
a_{_{k+1,1}} \\
u'&=a_{_{1,r+1}}\ldots a_{_{1,1}}b_{_{1,1}} \ldots
\overline{b_{_{k,k}}}a_{_{k,r+k}}\ldots a_{_{k,1}}b_{_{k,k-1}}\ldots
b_{_{k,1}}a_{_{k+1,r+k}}\ldots a_{_{k+1,1}} \end{aligned}
$$
where $a_{i,j}>0$ and $b_{i,j}<0$ above  and let  $S=P(u)$ and
$T=P(u')$.
\input{epsf}
\begin{figure}[h]
\vspace{-.4in}
\begin{center}
$\begin{array}{c} \epsfysize=7in \epsffile{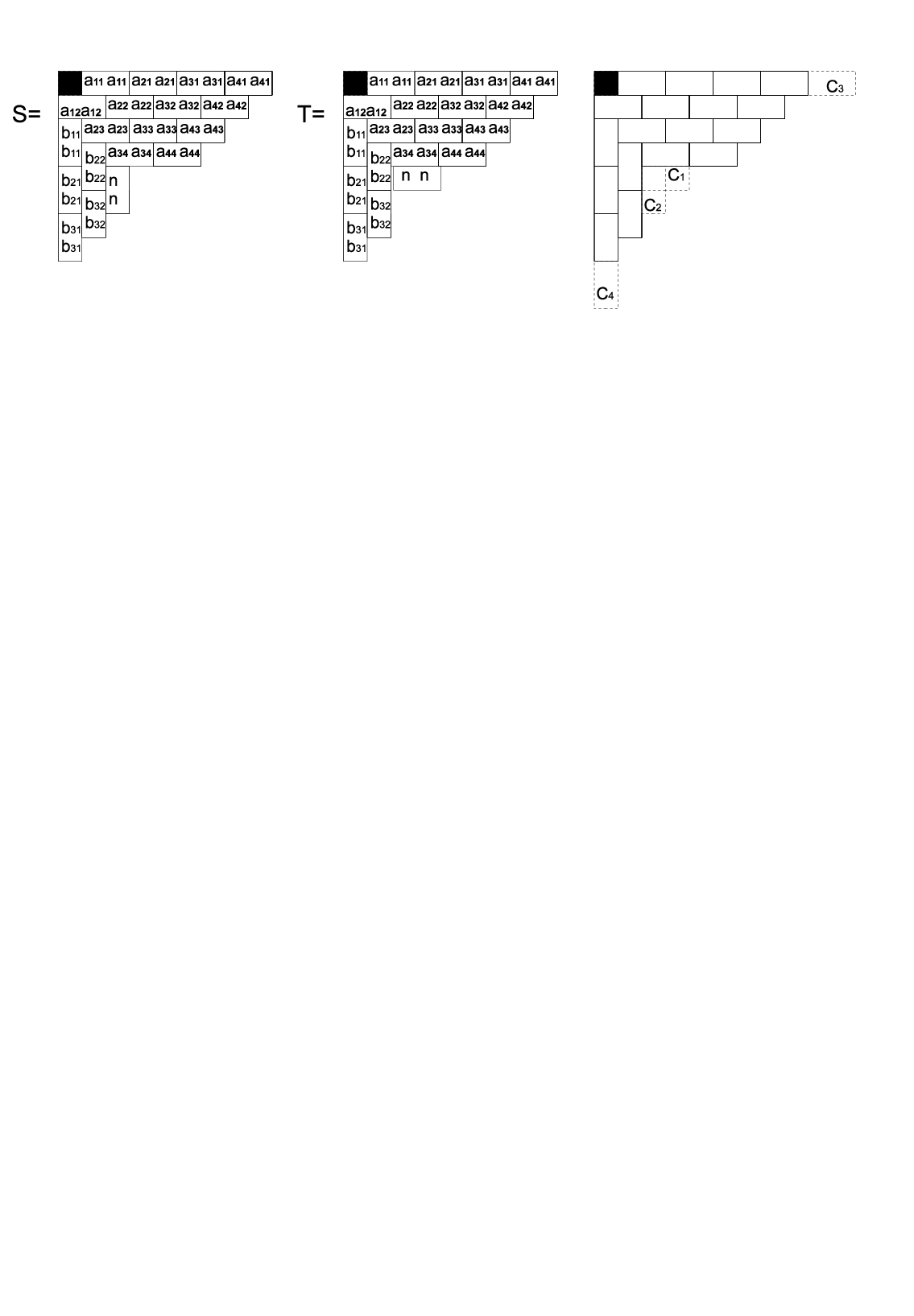}
\end{array}$
\end{center} \vspace{-5.7in}
\caption{ }\label{main.theorem.fig10}
\end{figure}

We first assume that $b_{_{k,k}}=\overline{n}$ and
$b_{_{k,1}}<z<a_{_{k+1,1}}$.   Then
 $S$ and $T$ differ only by the domino cell double labeled by $n$ i.e.
 $S_{<n}=T_{<n}$ with shape
 $$(s,s-1,\ldots,s-k, s-k-4,s-k-4,s-k-4,s-k-5,\ldots,2,1)$$
 for some $s=r+2k+2$. Therefore $S_{<n}=T_{<n}$ has exactly four empty domino
 corners:
$$C_1=\{(k+1,k),(k+1,k+1)\},~C_2=\{(k+1,k),(k+2,k)\},~C_3=\{(1,s+1),(1,s+2)\},~C_4= \{(s,1),(s+1,1)\}$$
as illustrated in Figure~\ref{main.theorem.fig10}   for $r=1$ and
$k=3$ (where the absolute value on $b_{i,j}$ is removed for the sake
of simplicity). Recall that
 by Lemma
~\ref{Garfinkle's.algorithm.lemma} one has
$$S^{\downarrow
z}=(S_{<n})^{\downarrow z} \leftarrow \mathrm{Dom}(S,n)~\mbox{ and
}~ T^{\downarrow z}=(T_{<n})^{\downarrow z} \leftarrow
\mathrm{Dom}(T,n)$$ where $\mathrm{Dom}(S,n)=[C_2,(n,n)]$ and
$[\mathrm{Dom}(T,n)=[C_1,(n,n)]$. On the other hand the assumption
$$b_{_{k,1}}<z<a_{_{k+1,1}}$$ yields   two choices for the new domino
cell appearing in $(S_{<n})^{\downarrow z}=(T_{<n})^{\downarrow z}$,
which  are $C_1$  and $C_2$. Now one can observe that whether $C_1$
or $C_2$ appears  sliding $\mathrm{Dom}(S, n)$ and $\mathrm{Dom}(T,
n)$ over $(S_{<n})^{\downarrow z}=(T_{<n})^{\downarrow z}$ gives the
same tableau. Therefore $P^r(uz)=P^r(u'z)$ as desired.

Now assume  $a_{_{k+1,r+k}}=n$. This time $S=P(u)$ and $T=P(u')$
differ only by the domino cell double labeled by $|b_{_{k,k}}|$, as
illustrated  in Figure~\ref{main.theorem.fig99} (a) and (b) for
$r=1$ and $k=3$.
\input{epsf}
\begin{figure}[h]
\vspace{-.4in}
\begin{center}
$\begin{array}{c} \epsfysize=7.2in \epsffile{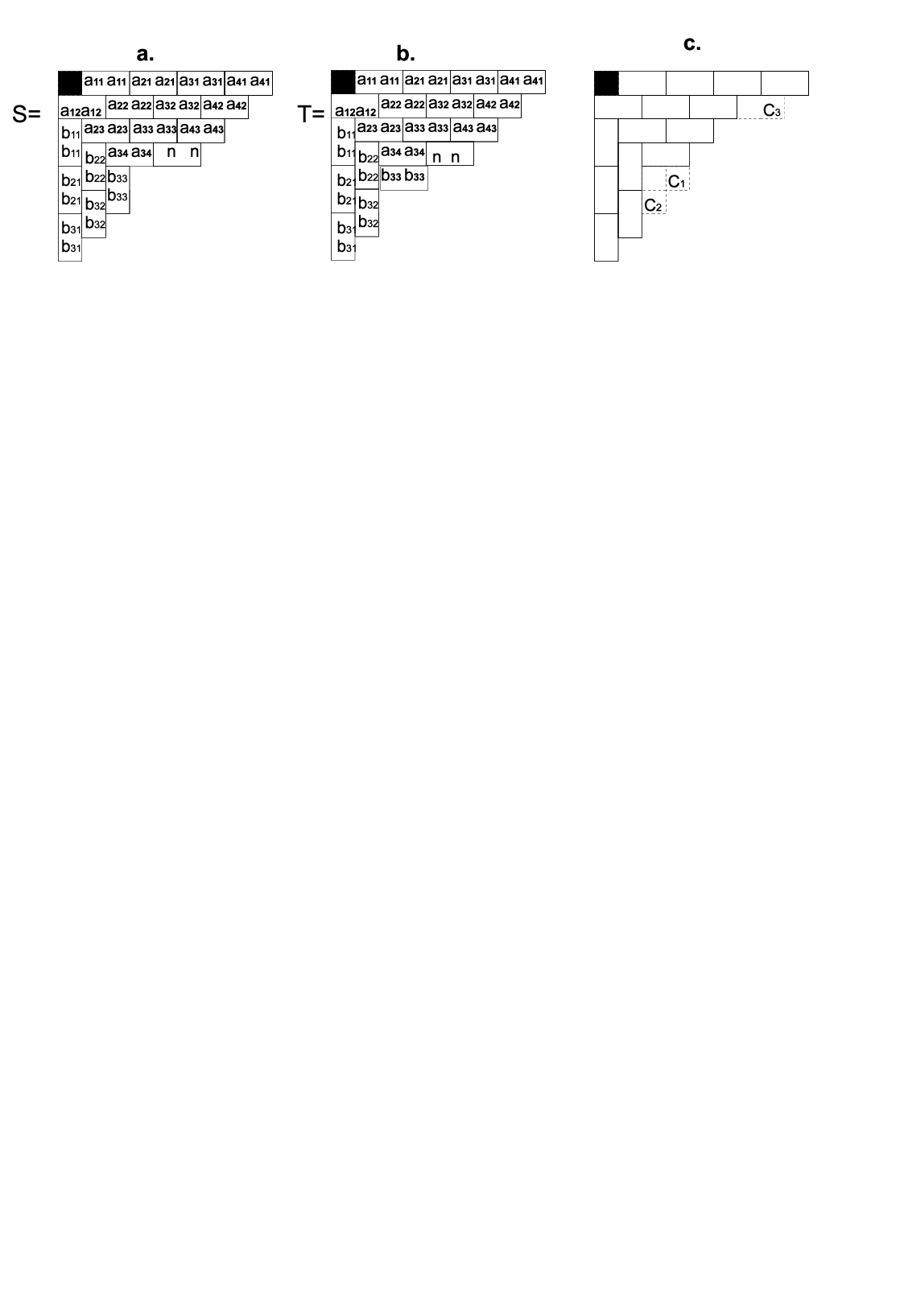}
\end{array}$
\end{center} \vspace{-5.9in}
\caption{ }\label{main.theorem.fig99}
\end{figure}

Let $\mathcal{A}$ be the set of labels which are greater
 than or equal to $|b_{_{k,k}}|$. Therefore $$\mathcal{A}=\{|b_{_{k,k}}|,a_{_{k+1,i}},\ldots,
a_{_{k+1,r+k-1}},a_{_{k+1,r+k}}=n\}$$ for some $1\leq i\leq r+k$.
Consider the tableaux $S'=S_{<|b_{_{k,k}}|}$ and
$T'=T_{<|b_{_{k,k}}|}$ which are obtained by erasing from $S$ and
$T$ respectively the domino cells double labeled by $\mathcal{A}$.
Now $S'=T'$ and moreover $(S')^{\downarrow z}=(T')^{\downarrow z}$.
Recall that  under the assumption $a_{_{k+1,r+k}}=n$ we have
$$ \begin{aligned}
&\mbox{ either}~~ b_{_{k,1}}<z<a_{_{k,1}}\\
&\mbox{ or}~~   ~a_{_{k,1}}<z<a_{_{k+1,1}}\mbox{ and}~~
a_{_{k+1,i}}<a_{_{k,i+1}} \mbox{ for some}~~1<i\leq k-1.
 \end{aligned}$$
Observe that in case $b_{_{k,1}}<z<a_{_{k,1}}$ we have two choices
for the new domino cell appearing in $(S')^{\downarrow
z}=(T')^{\downarrow z}$ which are, as illustrated in
Figure~\ref{main.theorem.fig99} (c),
$$C_1=\{(k+1,k),(k+1,k+1)\},~C_2=\{(k+1,k),(k+2,k)\}.$$
On the other hand if  $a_{_{k,1}}<z<a_{_{k+1,1}}$     one would add
the domino cell $C_3=\{(i,r+2k+1),(i,r+2k+2)\}$ in the above list
but existence of $C_3$  requires $a_{_{k,1}}<z<a_{_{k+1,1}}$ and
$a_{_{k+1,i}}\not<a_{_{k,i+1}}$ for any $1<i\leq k-1$. Therefore in
both cases $C_1$ and $C_2$ are the only choices for a new domino
cell appearing in $(S')^{\downarrow z}=(T')^{\downarrow z}$.

Recall that
$$\begin{aligned}
S^{\downarrow z}=&(S')^{\downarrow z}\leftarrow
\mathrm{Dom}(S,|b_{k,k}|) \leftarrow \mathrm{Dom}(S,a_{k+1,i})\ldots
\leftarrow
\mathrm{Dom}(S,a_{k+1,r+k-1}) \leftarrow \mathrm{Dom}(S,n)\\
T^{\downarrow z}=&(T')^{\downarrow z}\leftarrow
\mathrm{Dom}(T,|b_{k,k}|)\leftarrow \mathrm{Dom}(T,a_{k+1,i})\ldots
\leftarrow \mathrm{Dom}(T,a_{k+1,r+k-1})   \leftarrow
\mathrm{Dom}(T,n)
\end{aligned}$$
where  $\mathrm{Dom}(S,|b_{k,k}|)=C_2 $ and
$\mathrm{Dom}(T,|b_{k,k}|)=C_1$.  Now whatever $C_1$ or $C_2$
appears in  $(S')^{\downarrow z}=(T')^{\downarrow z}$,  sliding the
domino cell double labeled by $|b_{k,k}|$ gives the same tableau,
i.e.,
$$(S')^{\downarrow z}\leftarrow
\mathrm{Dom}(S,|b_{k,k}|)=(T')^{\downarrow z}\leftarrow
\mathrm{Dom}(T,|b_{k,k}|).$$ Moreover since
$\mathrm{Dom}(S,a_{k+1,j})= \mathrm{Dom}(T,a_{k+1,j})$  for all
$i\leq j\leq r+k$ we have $P^r(uz)=S^{\downarrow z}=T^{\downarrow
z}=P^r(u'z)$ as desired.
\end{proof}

\begin{theorem}\label{main.theorem2} If $\alpha$ and $\beta$
 have the same insertion $r$-domino tableaux then $\alpha~\rp~\beta$ in
 $B_n$.
 \end{theorem}

\begin{proof}  We will proceed by induction.  If $n=1$ there is nothing to prove, so suppose that
 the statement holds for all signed permutations of size $n-1$.

Let   $\al=\al_1\ldots\al_{n-1}\al_n$ and
$\beta=\beta_1\ldots\beta_{n-1}\beta_n$ satisfies
$T=P^r(\alpha)=P^r(\beta)$.  Therefore there exist two domino
corners say $A$ and $B$ of $T$ such that
\begin{equation}\label{main.theorem.eq1}
\begin{aligned}
T^{\uparrow A}=P^r(\alpha_1\ldots \alpha_{n-1}) ~\text{ and }~
\eta(T^{\uparrow A})=\alpha_n\\
T^{\uparrow B}=P^r(\beta_1\ldots\beta_{n-1})~\text{ and }~
\eta(T^{\uparrow B})=\beta_n. \end{aligned}
\end{equation}

In the following we suppose that $A$ lies below $B$.
\vskip.1in
 \noindent {\bf Cases 1:}  $A=B$. Then clearly  $T^{\uparrow A}=T^{\uparrow B}$
 and  $\alpha_n=\eta(T^{\uparrow A})=\eta(T^{\uparrow B})= \beta_n$.  Then since $$P(\alpha_1\ldots \alpha_{n-1} )=T^{\uparrow A}=T^{\uparrow B}=P(\beta_1\ldots\beta_{n-1})$$  we have, by induction, that
$\alpha_1\ldots \alpha_{n-1} ~\rp~ \beta_1\ldots\beta_{n-1}$ and
therefore $\alpha=\alpha_1\ldots \alpha_{n-1}\alpha_n
~\rp~\beta_1\ldots\beta_{n-1}\beta_n=\beta$  as desired.

\input{epsf}
\begin{figure}[h]
\vspace{-.4in}
\begin{center}
$\begin{array}{c} \epsfysize=8in \epsffile{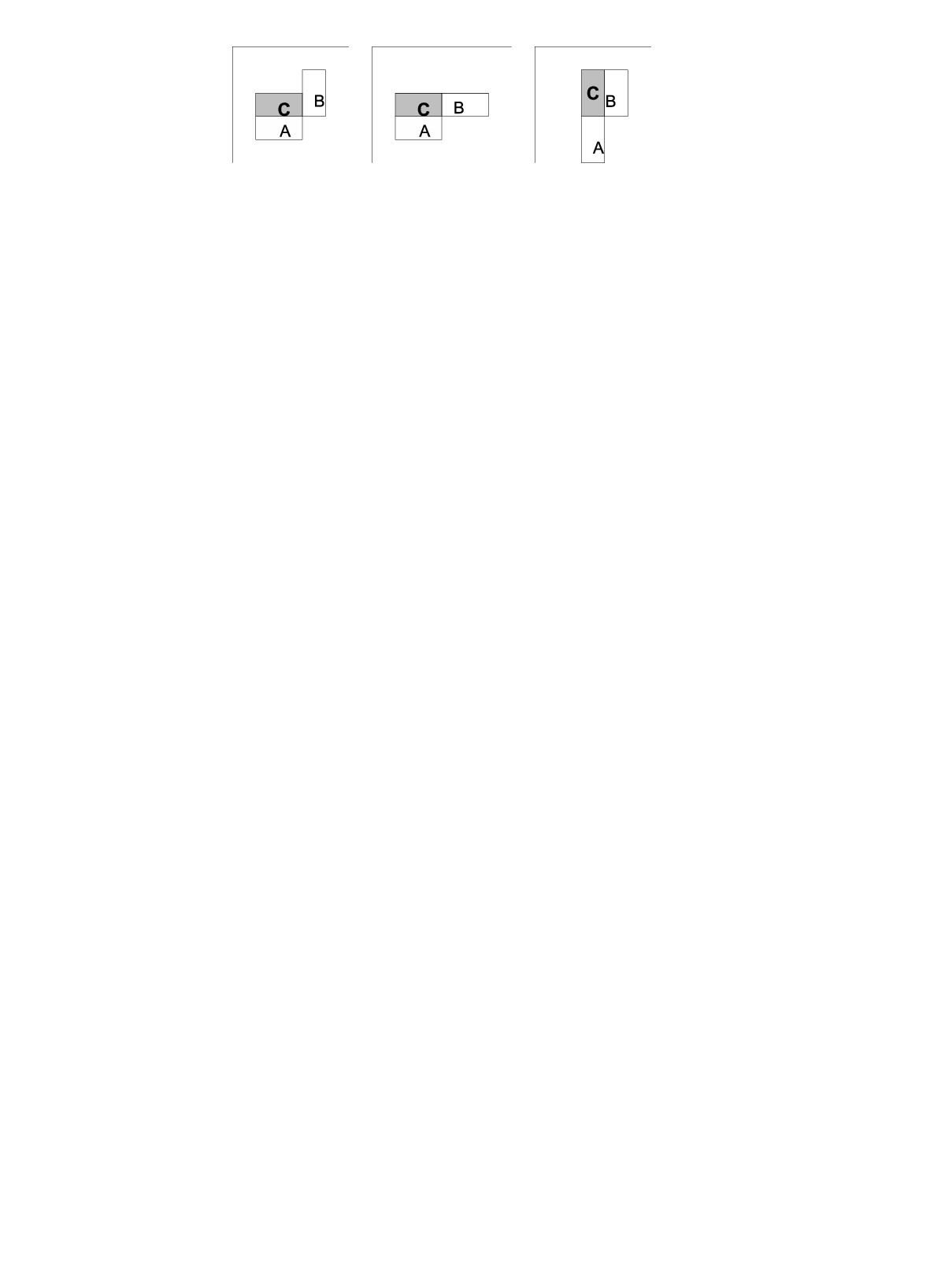}
\end{array}$
\end{center} \vspace{-7.2in}
\caption{Some illustrations  for Case 2}\label{main.theorem.fig0}
\end{figure}

\noindent {\bf Cases 2:} $A\not= B$ and $(T,A,
\mathrm{ne})\cap(T,B,\mathrm{sw})$   contains a domino corner, say
$C$ as illustrated in  Figure \ref{main.theorem.fig0}. Let
$$ \begin{aligned}
&b=\eta(T^{\uparrow A \uparrow B})\\
&c=\eta(T^{\uparrow A \uparrow B \uparrow C}) \end{aligned}$$
 and
let $\tilde{u}$ be a signed word such that
$P^r(\tilde{u})=T^{\uparrow A \uparrow B \uparrow C}$. Therefore
$$P^r(\tilde{u}cb\alpha_n)={P^r(\tilde{u})}^{\downarrow c\downarrow b\downarrow\alpha_n}=
(T^{\uparrow A \uparrow B \uparrow C})^{\downarrow c\downarrow
b\downarrow \alpha_n}=(T^{\uparrow A \uparrow B})^{\downarrow
b\downarrow \alpha_n}=(T^{\uparrow A})^{\downarrow \alpha_n}=T$$ and
by induction hypothesis $\tilde{u}cb~ \rp~ \alpha_1\ldots
\alpha_{n-1}$ since $P^r(\tilde{u}cb)= T^{\uparrow
A}=P^r(\alpha_1\ldots \alpha_{n-1})$ . Therefore
$$\tilde{u}cb\alpha_n ~\rp
~ \alpha.$$

Observe that since $P^r(\tilde{u})=T^{\uparrow A \uparrow B \uparrow
C}$,  the recording tableau $Q^r( \tilde{u}cb\alpha_n)$ has its
domino cells $A$, $B$ and $C$  labeled with $(n,n)$, $(n-1,n-1)$ and
$(n-2,n-2)$ respectively.

On the other hand  having  $B$  in $(T,A, \mathrm{ne})$ and $C$ in
$(T,B, \mathrm{sw})$ yields   by Lemma~\ref{reverse.insertion.lem2} that
$$b=\eta(T^{\uparrow A
\uparrow B})> \eta(T^{\uparrow A})= \alpha_n ~\text{ and
}~b=\eta(T^{\uparrow A \uparrow B})>\eta(T^{\uparrow A \uparrow B
\uparrow C})= c.$$ Therefore by Corollary~\ref{Barbash.Vogan.cor2},
$$\begin{aligned}
\text{either}& ~ b>\alpha_n>c, ~\text{ and hence }~
u=\tilde{u}cb\alpha_n ~\rda~ \tilde{u}bc\alpha_n=w ~\text{ and }
~V_{n-1,n-2}(Q^r(u))=Q^r(w) \\
\text{or}&~ b>c>\alpha_n, ~\text{ and hence }~ u=\tilde{u}cb\alpha_n
~\rda~ \tilde{u}c\alpha_nb=w~\text{ and } ~
V_{n-1,n-2}(Q^r(u))=Q^r(w)
\end{aligned}
$$
The last argument  implies that in both cases the signed permutation
$w$ has its recording tableau $Q^r(w)$ obtained by interchanging the
labels $(n,n)$ of $A$ and $(n-1,n-1)$ of $B$ in $Q^r(u)$ which means
 that $Q^r(w)$ has its domino corner $B$ labeled with $(n,n)$.
 So we have
$$P^r(w_1\ldots
w_{n-1})=T^{\uparrow B}=P^r(\beta_1\ldots\beta_{n-1}) ~\text{ and
}w_n=\beta_n. $$
 Now by induction, $w_1\ldots w_{n-1}
~\rp~\beta_1\ldots\beta_{n-1}$ and  therefore  $w ~\rp~ \beta$.
Hence  $\alpha~\rp~u~\rp~w~\rp~\beta$ as required.

\vskip.1in \noindent {\bf Case 3:} $A\not=B$,  $(T,A,
\mathrm{ne})\cap(T,B,\mathrm{sw})$ is a staircase shape
$(s,s-1,\ldots,1)$ for $s\geq1$ and  $A\cap B$ is a single box. The
condition $A\cap B$ is a single box forces that $s=1$. There are
several subcases.

\vskip.1in \noindent {\bf Case 3.1:}   We assume that $T$ has no
domino corner beyond $A$ and $B$. Let $\lambda$ and $\lambda'$ be
respectively  the  smallest and the  largest rectangular shape
containing both $A$ and $B$ whose east and south boundary coincides
with the boundary of $\sh(T)$. Then clearly $\lambda=(2,2)$ and we
have either $\lambda\subsetneq \lambda'$ or $\lambda= \lambda'$ as
illustrated in Figure~\ref{main.theorem.fig21}.
\input{epsf}
\begin{figure}[h]
\vspace{-.2in}
\begin{center}
$\begin{array}{c} \epsfysize=5.4in \epsffile{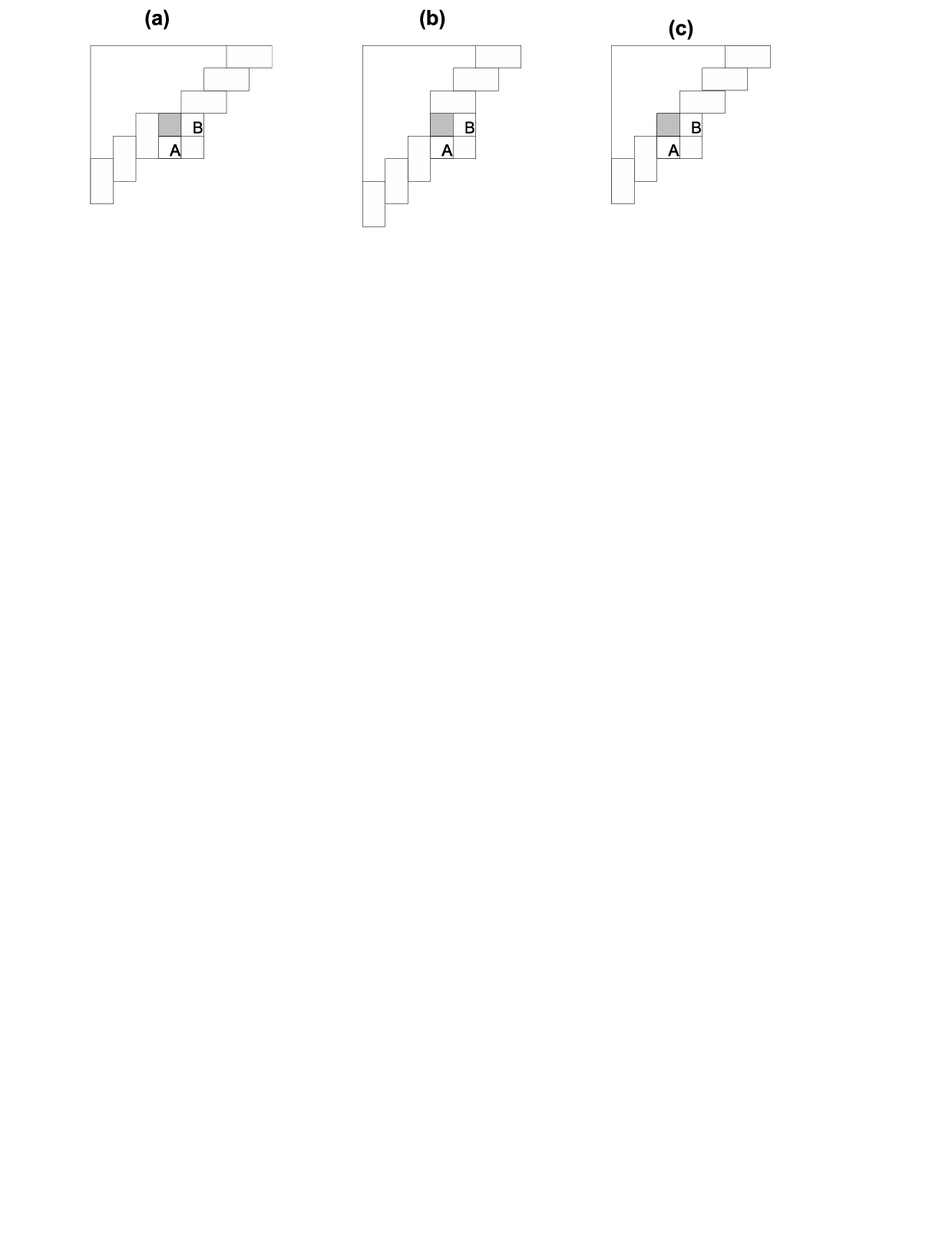}
\end{array}$
\end{center} \vspace{-4.6in}
\caption{Case 3.1: $T$ has no domino corner beyond $A$ and $B$ and
the partitions  $(3,3)$, $(2,2,2)$ and $(2,2)$ determine $\lambda'$
in (a),(b) and (c) respectively.}\label{main.theorem.fig21}
\end{figure}
\vskip.1in \noindent {\bf Case 3.1.1:} We first suppose that
$\lambda\subsetneq \lambda'$ as  illustrated in
Figure~\ref{main.theorem.fig21}($a$) and ($b$). Since the other case
can be dealt with in the same manner after taking the transpose of
$T$, below we just consider the case
Figure~\ref{main.theorem.fig21}($a$), where there exist a vertical
domino cell to the left of $\lambda$ in $\lambda'$.

\input{epsf}
\begin{figure}[h]
\vspace{-1.5in}
\begin{center}
$\begin{array}{c} \epsfysize=7in \epsffile{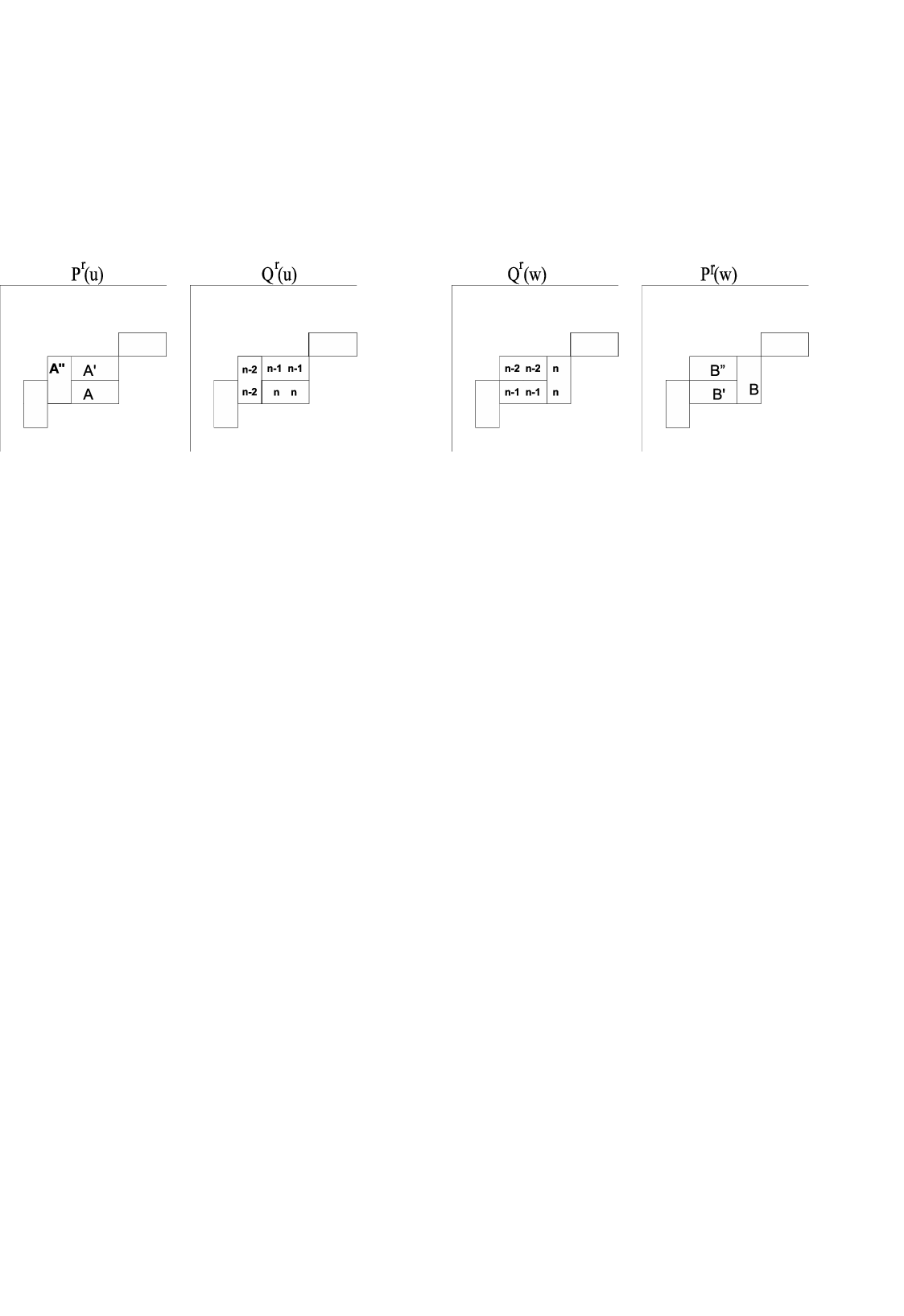}
\end{array}$
\end{center} \vspace{-4.8in}
\caption{Case 3.1.1}\label{main.theorem.fig3}
\end{figure}

Now observe through Figure~\ref{main.theorem.fig3} that   we have a
domino corner $A'$ of $T^{\uparrow A}$ and $A''$ of $T^{\uparrow A
\uparrow A'}$ as given in Figure~\ref{main.theorem.fig3}. Let
$a'=\eta(T^{\uparrow A \uparrow A'})~\mbox{ and }~
a''=\eta(T^{\uparrow A\uparrow A'\uparrow A''}).$ Suppose
$\tilde{u}$ be a signed word such that $P^r(\tilde{u})=T^{\uparrow
A\uparrow A'\uparrow A''}$. Then the signed permutation $u=
\tilde{u}a''a'\alpha_n$ has $P^r(u)=T$ whereas  its recording
tableau  $Q^r(u)$ must have the form  as it is shown in
Figure~\ref{main.theorem.fig3}.

Furthermore  having  $A'$  in $(T,A, \mathrm{ne})$ and $A''$ in
$(T,A', \mathrm{sw})$ yields   by Lemma~\ref{reverse.insertion.lem2}
that
$$a'=\eta(T^{\uparrow A
\uparrow A'})> \eta(T^{\uparrow A})= \alpha_n ~\text{ and
}~a'=\eta(T^{\uparrow A \uparrow A'})>\eta(T^{\uparrow A \uparrow A'
\uparrow A''})= a''.$$ Therefore we have
$$\begin{aligned}
\text{either}& ~ a''<\alpha_n<a', ~\text{ and hence }~
u=\tilde{u}a''a'\alpha_n ~\rda~ \tilde{u}a'a''\alpha_n=w ~\text{ and}~Q^r(w)= V_{n-2,n-1}(Q^r(u)) \\
\text{or}&~ \alpha_n<a''<a', ~\text{ and hence
}~u=\tilde{u}a''a'\alpha_n ~\rda~ \tilde{u}a''\alpha_na'=w~\text{
and } ~Q^r(w)= V_{n-2,n-1}(Q^r(u)).
\end{aligned}
$$
In both cases Corollary~\ref{Barbash.Vogan.cor2} yields that the
recording tableau $Q^r(w)$ of $w$ has the form  in
Figure~\ref{main.theorem.fig3}.

Now  since $P^r( \tilde{u}a''a')=T^{\uparrow
A}=P^r(\alpha_1\ldots\alpha_{n-1})$, we have  by induction
$\tilde{u}a''a'~\rp~\alpha_1\ldots\alpha_{n-1}$. Therefore
$$u=\tilde{u}a''a'\alpha_n~\rp~\alpha_1\ldots\alpha_{n-1}\alpha_n=\alpha.$$
Similarly since $P^r(w_1\ldots w_{n-1})=T^{\uparrow
B}=P^r(\beta_1\ldots\beta_{n-1})$, by induction  $w_1\ldots
w_{n-1}~\rp~\beta_1\ldots\beta_{n-1}$. On the other hand since
$w_n=\beta_n$ we have
 $w_1\ldots
w_{n-1}\beta_n~\rp~\beta_1\ldots\beta_{n-1}\beta_n=\beta$. Hence
 $$\alpha~\rp~u~\rda~w~\rp~\beta.$$

\vskip.1in \noindent {\bf Case 3.1.2:} No we suppose that $\lambda=
\lambda'$ as   illustrated in Figure~\ref{main.theorem.fig21}($c$).
Observe through Figure~\ref{main.theorem.fig4}  that the grey area
in the first tableau  has a staircase shape and  since there are no
other domino corners of $T$, we must have either $A$ or $B$ labeled
by $(n,n)$.
\input{epsf}
\begin{figure}[h]
\vspace{-.4in}
\begin{center}
$\begin{array}{c} \epsfysize=5.5in \epsffile{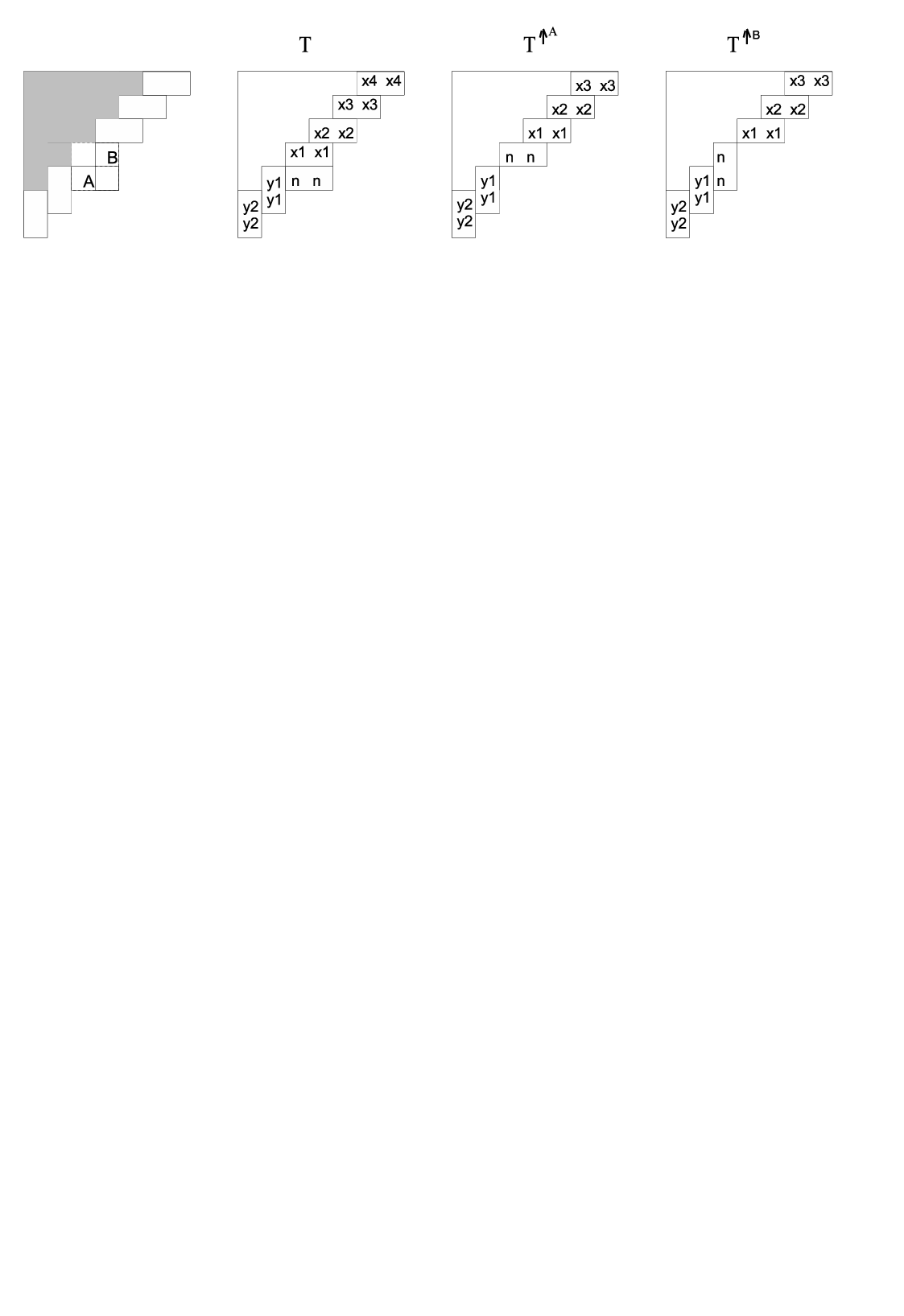}
\end{array}$
\end{center} \vspace{-4.6in}
\caption{Case 3.1.2 }\label{main.theorem.fig4}
\end{figure}

Suppose that the horizontal domino cell  $A$ is labeled by $(n,n)$
(The other case can be also dealt with taking the transpose of the
tableau). So as Figure~\ref{main.theorem.fig4} illustrates, we have
only the horizontal domino cells double labeled by $n, x_1, \ldots,
x_k$ and  the vertical domino cells double labeled by $y_1,\ldots,
y_l$ where $k\geq 0$, $l\geq 0$ and   $$k+l=r+1.$$ Therefore at
least one of $k$ and $l$ must be nonzero.  W.L.O.G. assume that
$k\geq 1$.  In this case observe that $\eta(T^{\uparrow
A})=\eta(T^{\uparrow B})=x_k>0$ and this yields
$\alpha_n=\beta_n=x_k$. Let $\tilde{u}= \overline{y_1}\ldots
\overline{y_l}x_1\ldots x_{k-1}$. Clearly
$P^r(n\tilde{u}x_k)=T=P^r(\overline{n}\tilde{u}x_k)$ and
$$n\tilde{u}x_k~\rdc~\overline{n}\tilde{u}x_k.$$
On the other hand
$P^r(n\tilde{u})=T^{\uparrow A}$ and $P^r(\overline{n} \tilde{u})=T^{\uparrow
B}$ and  by induction hypothesis we have
$n\tilde{u}~\rp~\alpha_1\ldots\alpha_{n-1}$ and
$\overline{n}\tilde{u}~\rp~\beta_1\ldots\beta_{n-1}$. Hence
$$\alpha=\alpha_1\ldots\alpha_{n-1}x_k ~\rp~ n\tilde{u}x_k~
\rdc~\overline{n}\tilde{u}x_k ~\rp~\beta_1\ldots\beta_{n-1}x_k=\beta.
$$

\noindent {\bf Case 3.2:}  Suppose that $T$ has another domino
corner, say $C$. Then $C$ must lie either  in $(T,B, \mathrm{ne})$
or $(T,A, \mathrm{sw})$. Below we assume that  $C$  lies   in $(T,B,
\mathrm{ne})$, since the other case can be dealt with in the same
manner by  considering the transpose of $T$.

\noindent {\bf Case 3.2.1:} We first suppose that
 $(T,C,\mathrm{sw})\cap(T,B,\mathrm{ne})$ contains a domino corner as illustrated in
 Figure~\ref{main.theorem.fig22}.
Let $\sigma $ be a permutation such that $ T^{\uparrow
C}=P^r(\sigma_1\ldots \sigma_{n-1})$ and $\eta(T^{\uparrow
C})=\sigma_n$. Since also $(\mathrm{T,C,sw})\cap(T,A,\mathrm{ne})$
contains a domino corner,  $A$  and $C$ satisfy Case 2 and this
gives  $\alpha~\rp~ \sigma$ and $\beta~\rp~ \sigma$. Therefore
$\beta ~\rp~ \alpha$.

\input{epsf}
\begin{figure}[h]
\vspace{-.3in}
\begin{center}
$\begin{array}{c} \epsfysize=5.7in \epsffile{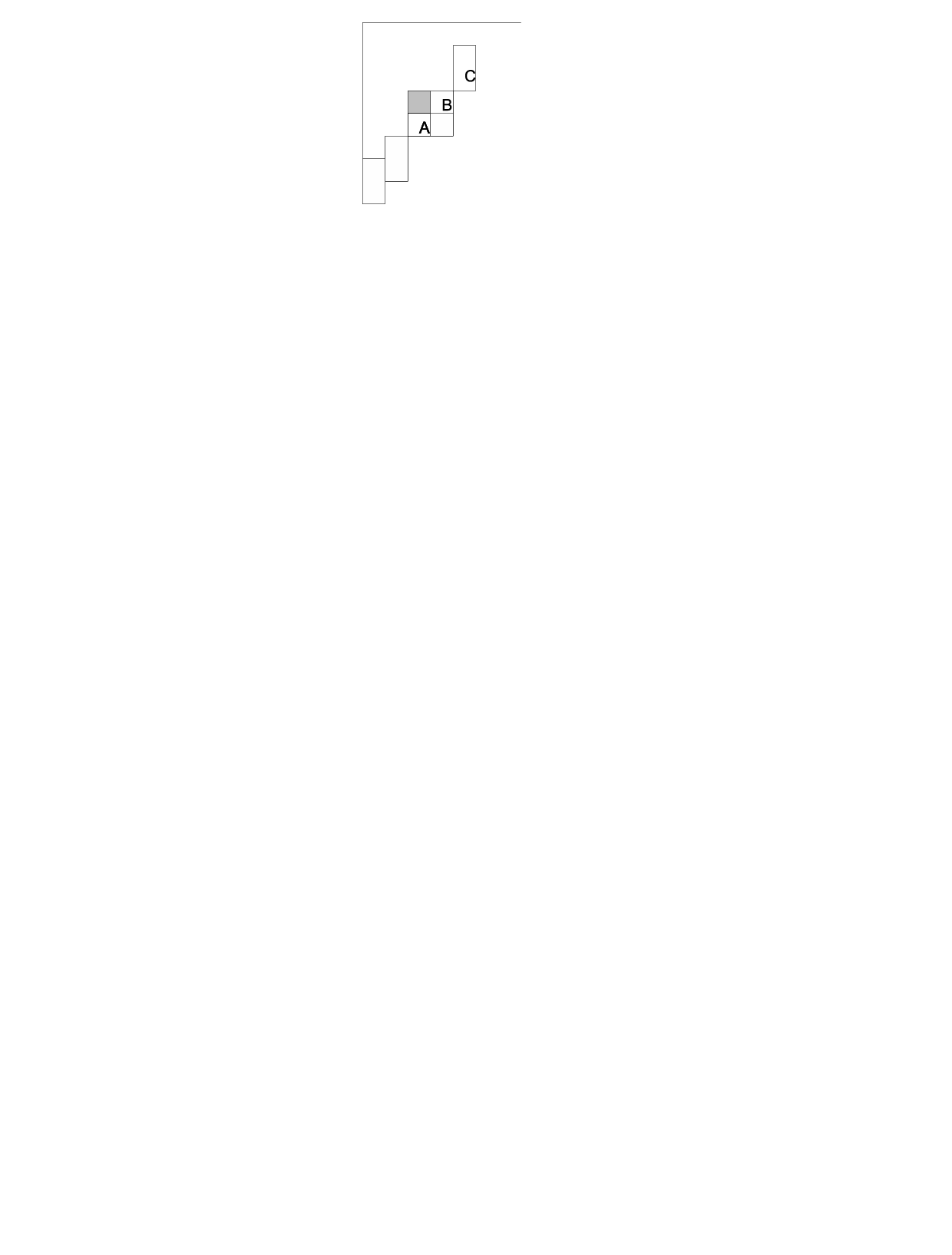}
\end{array}$
\end{center} \vspace{-5in}
\caption{Case 3.2.1 }\label{main.theorem.fig22}
\end{figure}

\vskip.1in \noindent {\bf Case 3.2.2:} Now we  suppose  that  $C$
satisfies
$$(T,C,\mathrm{sw})\cap(T,B,\mathrm{ne})=(s,s-1,\ldots,2,1)~\mbox{ for some } s\geq 1.$$
Further we can assume that $C$ is the only domino corner lying in
$(T,B,\mathrm{ne})$, since  any other domino corner in
$(T,B,\mathrm{ne})$ must satisfy  Case 3.2.1.  Below,
Figure~\ref{main.theorem.fig23} illustrates the possible subcases:
\input{epsf}
\begin{figure}[h]
\vspace{-.2in}
\begin{center}
$\begin{array}{c} \epsfysize=6in \epsffile{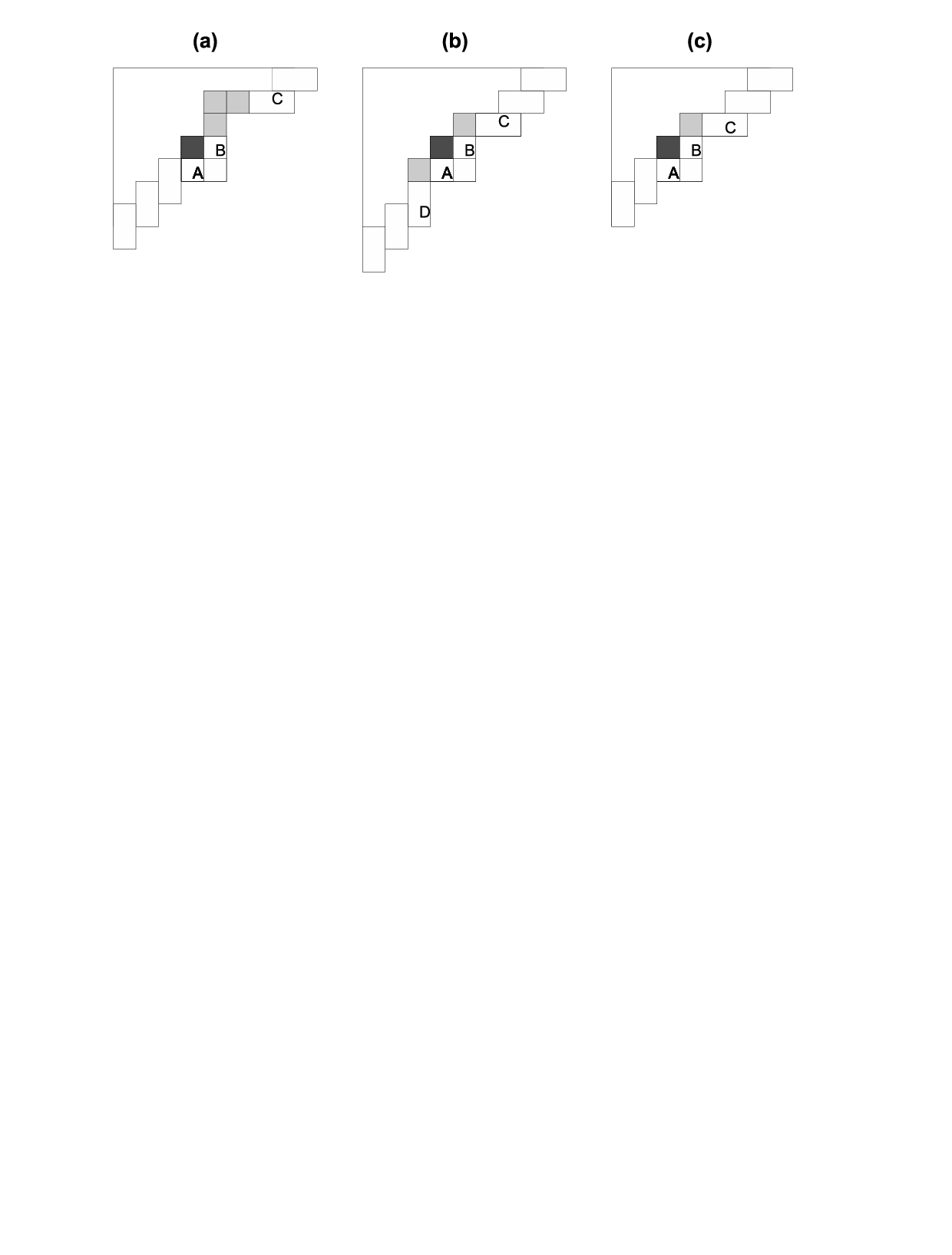}
\end{array}$
\end{center} \vspace{-4.9in}
\caption{Case 3.2.2 }\label{main.theorem.fig23}
\end{figure}

We first suppose  that  $(T,C,\mathrm{sw})\cap(T,B,\mathrm{ne})$ is
a staircase shape $(s,s-1,\ldots,1)$ for $s>1$ as illustrated in
Figure~\ref{main.theorem.fig23}($a$). Then whether
$(T,A,\mathrm{sw})$
 contains another domino corner   or not, $A$ and $B$ are contained
in a rectangular shape which is strictly larger than
$\lambda=(2,2)$. Therefore   this case is similar to the one
 studied in Case 3.1.1.,  and the same method applied there gives  $\alpha~\rp~ \beta$.

Now  suppose that   $(T,C,\mathrm{sw})\cap(T,B,\mathrm{ne})$ is a
single box and that there exist a domino corner, say $D$ in
$(T,A,\mathrm{sw})$. If $(T,D,\mathrm{ne})\cap(T,A,\mathrm{sw})$
contains a domino corner then one can apply the same method of Case
3.2.1 for the domino corners $A$ and $D$ to get the desired result.
On the other hand the case that
$(T,D,\mathrm{ne})\cap(T,A,\mathrm{sw})$ contains a staircase shape
$(s,s-1,\ldots,1)$ for some $s>1$,  is similar to the one pictured
in Figure~\ref{main.theorem.fig23}($a$), therefore $\alpha~\rp~
\beta$ follows directly. Now the case that
$(T,D,\mathrm{ne})\cap(T,A,\mathrm{sw})$ is a single box is
illustrated in Figure~\ref{main.theorem.fig23}($b$). For this case
let $\sigma$ and $\delta$ be two permutations in $B_n$ satisfying:
$$
T^{\uparrow C}=P^r(\sigma_1\ldots \sigma_{n-1}),~ \eta(T^{\uparrow
C})=\sigma_n~\text{ and }~T^{\uparrow D}=P^r(\delta_1\ldots
\delta_{n-1}), ~ \eta(T^{\uparrow D})=\delta_n.
$$
Observe that both $(T,C,\mathrm{sw})\cap(T,A,\mathrm{ne})$ and
$(T,D,\mathrm{ne})\cap(T,B,\mathrm{sw})$ contain a domino corner
i.e., the pair $A$ and $C$ and  similarly the pair $B$ and $D$
satisfy Case 2. Therefore   $\alpha~\rp~ \sigma$ and $\beta~\rp~
\delta$. On the other hand $C$ and $D$ also satisfy Case 2, so  we have
$\sigma~\rp~ \delta$.  Hence the result  $\alpha~\rp~ \beta$
follows.

Lastly we suppose that $(T,C,\mathrm{sw})\cap(T,B,\mathrm{ne})$ is a
single box and that there exist no  domino corners in
$(T,A,\mathrm{sw})$ as  illustrated in
Figure~\ref{main.theorem.fig23}($c$) (See also
Figure~\ref{main.theorem.fig13} for  possible other variations,
including the transpose of $T$). Therefore
 $T$ has
shape
$$(s+i,s+i-1,\ldots,s+1,\underline{s,s-2,s-2,s-4},s-3,\ldots,2,1) ~\mbox{for some }~ s\geq 4,
~i\geq 0
$$
and it has  three domino corners
$$A=\{(i+3,s-3),(i+3,s-2)\},B=\{(i+2,s-2),
(i+3,s-2)\},C=\{(i+1,s-1),(i+1,s)\}
$$
as Figure~\ref{main.theorem.fig12}($a$) illustrates.
\begin{figure}[h]
\vspace{-.4in}
\begin{center}
$\begin{array}{c} \epsfysize=7in \epsffile{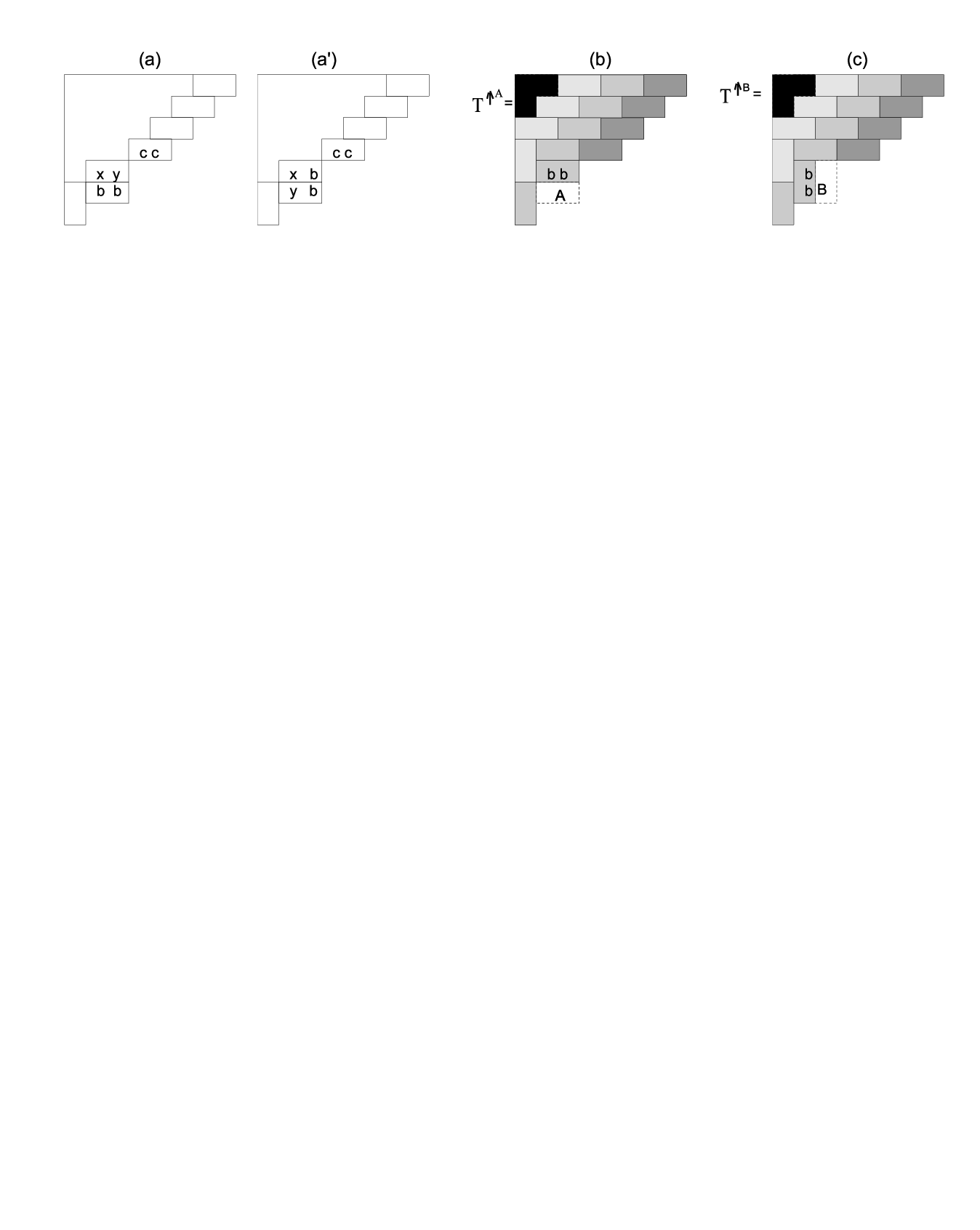}
\end{array}$
\end{center} \vspace{-5.9in}
\caption{Case 3.2.2, Figure~\ref{main.theorem.fig23}
(c)}\label{main.theorem.fig12}
\end{figure}

First observe that  $C$ must be double labeled by some number $c$,
since otherwise  none of the right most horizontal domino cells
above $C$ can be double labeled and  this contradicts to the fact
that $T$ is a standard $r$-domino tableau.

 Now we suppose that $A$
is double labeled by some  number $b$ and consider the horizontal
domino cell  $D=\{(i+2,s-3),(i+2,s-2)\}$ which is indicated by the
labeling $(x,y)$ in Figure~\ref{main.theorem.fig12}($a$). Observe
from Corollary~\ref{reverse.insertion.cor1} that Garfinkle's reverse
insertion algorithm applied on $A$ first pushes back $[D,(x,y)]$ and
then labels $D$ by $(b,b)$ as Figure~\ref{main.theorem.fig12}($b$)
suggests. On the other hand   Garfinkle's reverse insertion
algorithm applied on $B$ first pushes back $[D,(x,y)]$ and then
labels the vertical domino cell $\{(i+2,s-3),(i+3,s-3)\}$ by $(b,b)$
as Figure~\ref{main.theorem.fig12}($c$) illustrates. Since the
domino cell $D$ of $T$ is to be pushed back in the first step of
both reverse insertion, we have
$$\alpha_n=\eta(T^{\uparrow A})=z=\eta(T^{\uparrow B})=\beta_n
$$
and furthermore the resulting tableaux $T^{\uparrow A}$ and
$T^{\uparrow B}$ just differ by their domino cells labeled by
$(b,b)$. In fact one can get the same result if $B$ is double
labeled by a number $b$, just that this time the vertical domino
cell $D'=\{(i+2,s-3),(i+3,s-3)\}$ (that
Figure~\ref{main.theorem.fig12}($a'$) indicates by the labeling
$(x,y)$) is to be pushed back  in the first step of Garfinkle's
reverse insertion applied on $A$ and respectively $B$.

Observe that  the east most horizontal  and the south most vertical
domino cells of $T^{\uparrow B}$ must be double labeled. We first
apply the reverse insertion on the east most horizontal domino cells
starting from to bottom cell $C$ to top so that a sequence of
positive increasing numbers is obtained.  Next the reverse insertion
applied on the south most domino cells from right to left gives a
sequence of increasing negative numbers. Continuing the reverse
insertion of the east most horizontal domino cells and the south
most vertical domino cells in the remaining tableaux
 one at a time,  we end up by a
staircase shape $(r,r-1,\ldots,0)$ for some $r\geq 0$.

If the shape $(r,r-1,\ldots,0)$ is obtained by reverse inserting the
east most horizontal domino cells at the end, as Figure~\ref
{main.theorem.fig12} and Figure~\ref{main.theorem.fig13}(a)
illustrate, then last sequence obtained in this manner must be also
positive decreasing as the first sequence. Therefore for some $k\leq
1$ satisfying $s+i=r+2(k+1)$ we have the following word
$$u=a_{1,r+1}\ldots a_{1,1}b_{1,1}\ldots
a_{k,r+k}\ldots a_{k,1}b_{k,k}\ldots b_{k,1}a_{k+1,r+k}\ldots
a_{k+1,1} $$ where $a_{i,r+i}\ldots a_{1,1}$  represents positive
decreasing sequence obtained by reverse inserting east most
horizontal domino cells (Observe that $i=1$ represent the last
sequence)  and $b_{i,i}\ldots b_{1,1}$ represent negative increasing
sequence obtained by reverse inserting south most vertical domino
cells. Therefore
 the numbers $a_{_{i,j}}$ and
$b_{_{i,j}}$ in $u$ satisfy the following conditions.
$$\begin{aligned}
&a_{i,j}>0, ~b_{i,j}< 0 \\
&a_{i,j-1}<a_{i,j}<a_{i+1,j} ~\mbox{ and }~  |b_{i,j-1}|<|b_{i,j}|<|b_{i+1,j}|\\
&|b_{i,i}|<a_{i+1,r+i+1}<|b_{i+1,i+1}|~\mbox{  for all }~
i=1,\ldots, k-1.
\end{aligned}$$
 Consider the following word which is obtained by taking
 $|b_{k,k}|$  (Observe that  $b_{k,k}=b$ in Figure~\ref{main.theorem.fig13})
 in front of $a_{k,r+k}$ in $u$.
$$u'=a_{1,r+1}\ldots a_{1,1}b_{1,1}\ldots
\overline{b_{k,k}}~a_{k,r+k}\ldots a_{k,1}b_{k,k-1}\ldots b_{k,1}
a_{k+1,r+k}\ldots a_{k+1,1}.$$ Therefore $P^r(u')=T^{\uparrow A}$
whereas $P^r(u)=T^{\uparrow B}$ and since $(T^{\uparrow
B})^{\downarrow z}=T=(T^{\uparrow A})^{\downarrow z}$ we have
$$P^r(uz)=T=P^r(u'z).
$$

 Now we have the following
analysis on $z=\eta(T^{\uparrow B})=\eta(T^{\uparrow A})$: Observe
that either $b$ or $c$ must be equal to $n$. If $b=n$ then  any
number $z$  between  $b_{k,1}$ and $a_{k+1,1}$ satisfies
$$(T^{\uparrow B})^{\downarrow z}=T=(T^{\uparrow A})^{\downarrow z}$$
where $|b_{k,1}|$ and $a_{k+1,1}$ respectively are the labels of
south most vertical  and right most horizontal  domino cells in both
$T^{\uparrow B}$ and  $T^{\uparrow A}$. If  $c=n$ then insertion of
the number $z$ in both $T^{\uparrow B}$ and  $T^{\uparrow A}$ can
not bump the domino cell $C$ which is double labeled by $n$ in to
the next row since then resulting tableau is not equal to  $T$.
Therefore in this case   either  $b_{k,1}<z<a_{k,1}$ or
$a_{k,1}<z<a_{k+1,1}$ and  $a_{k,i+1}>a_{k+1,i}<$ for some $1<i\leq
k-1$. As a result  two words  $uz$ and $u'z$ satisfy $D_4^r$
relation with all $a_{_{i,j}}>0$ and $b_{_{i,j}}<0$.

If the shape $(r,r-1,\ldots,0)$ is obtained by reverse inserting the
south  most vertical  domino cells at the end, as
Figure~\ref{main.theorem.fig13}(b) illustrates then  last sequence
obtained in this manner  must be negative  increasing  as opposed to
the first sequence. Therefore for some $k\leq 1$ satisfying
$s+i=r+2(k+1)+1$, we have the following word
$$u=a_{_{1,r+1}}\ldots
a_{_{1,1}}b_{_{1,1}}\ldots a_{_{k,r+k}}\ldots a_{_{k,1}}
b_{_{k,k}}\ldots b_{_{k,1}}a_{_{k+1,r+k+1}}\ldots
a_{_{k+1,1}}b_{_{k+1,k}}\ldots b_{_{k+1,1}} $$ where
$a_{1,r+1}\ldots a_{1,1}$ and $b_{_{k+1,k}}\ldots b_{_{k+1,1}}$
represent respectively the last (negative increasing) and the first
(positive decreasing) sequences obtained in this manner. Moreover
the numbers $a_{_{i,j}}$ and $b_{_{i,j}}$ in $u$ satisfy the
following conditions.
$$\begin{aligned}
&a_{_{i,j}}<0 ~\mbox { and }~b_{_{i,j}}>0 \\
&|a_{_{i,j-1}}|<|a_{_{i,j}}|<|a_{_{i+1,j}}| ~\mbox{ and }~  b_{_{i,j-1}}<b_{_{i,j}}<b_{_{i+1,j}}\\
& |a_{_{i,r+i}}|<b_{_{i,i}}<|a_{_{i+1,r+i+1}}|~\mbox{for all
}~i=1,\ldots, k.
\end{aligned}$$
Now one can easily check the following word
$$u'=a_{_{1,r+1}}\ldots a_{_{1,1}}b_{_{1,1}}\ldots a_{_{k,r+k}}\ldots a_{_{k,1}}
\overline{a_{_{k+1,r+k+1}}}~b_{_{k,k}}\ldots
b_{_{k,1}}a_{_{k+1,r+k}}\ldots a_{_{k+1,1}} b_{_{k+1,k}}\ldots
b_{_{k+1,1}}$$ satisfies   $P^r(u')=T^{\uparrow A}$ whereas
$P^r(u)=T^{\uparrow B}$. Moreover a similar analysis  on the number
$z$ shows that $z$ satisfies one of the hypothesis of $D_5^r$,
therefore  two words  $uz$ and $u'z$ satisfy $D_5^r$ relation with
all $a_{_{i,j}}<0$ and $b_{_{i,j}}>0$.

Now recall that $P(u)=T^{\uparrow A}=P(\alpha_1\ldots \alpha_{n-1})$
and $P(u')=T^{\uparrow B}=P(\beta_1\ldots \beta_{n-1})$.  So we have
$$\alpha_n=z=\beta_n$$  and moreover $u~\rp~\alpha_1\ldots
\alpha_{n-1}$ and $u'~\rp~\beta_1\ldots \beta_{n-1}$ by induction.
Therefore $\alpha ~\rp~ uz ~\rp~ u'z ~\rp ~\beta$ as desired.

Note that for the tableaux that  Figure~\ref{main.theorem.fig13}(c)
and (d) illustrates, we first  apply  the reverse insertion the
south most vertical domino cells starting from $C$ to the left and
in that case a sequence of negative decreasing numbers is obtained.
Moreover, according to the sign of the last sequence obtained in the
same manner  one get either  $D_4^r$ relation with all
$a_{_{i,j}}<0$ and $b_{_{i,j}}>0$
(Figure~\ref{main.theorem.fig13}(c)) or  $D_5^r$ relation with all
$a_{_{i,j}}>0$ and $b_{_{i,j}}<0$
(Figure~\ref{main.theorem.fig13}(d)).

\vskip.1in \noindent {\bf Case 4:} $A\not=B$,  $(T,A,
\mathrm{ne})\cap(T,B,\mathrm{sw})$ is a staircase shape
$(s,s-1,\ldots,1)$ for $s\geq1$ and  $A\cap B$ is empty.
Figure~\ref{main.theorem.fig5} shows  several subcases.
\input{epsf}
\begin{figure}[h]
\vspace{-.3in}
\begin{center}
$\begin{array}{c} \epsfysize=6in \epsffile{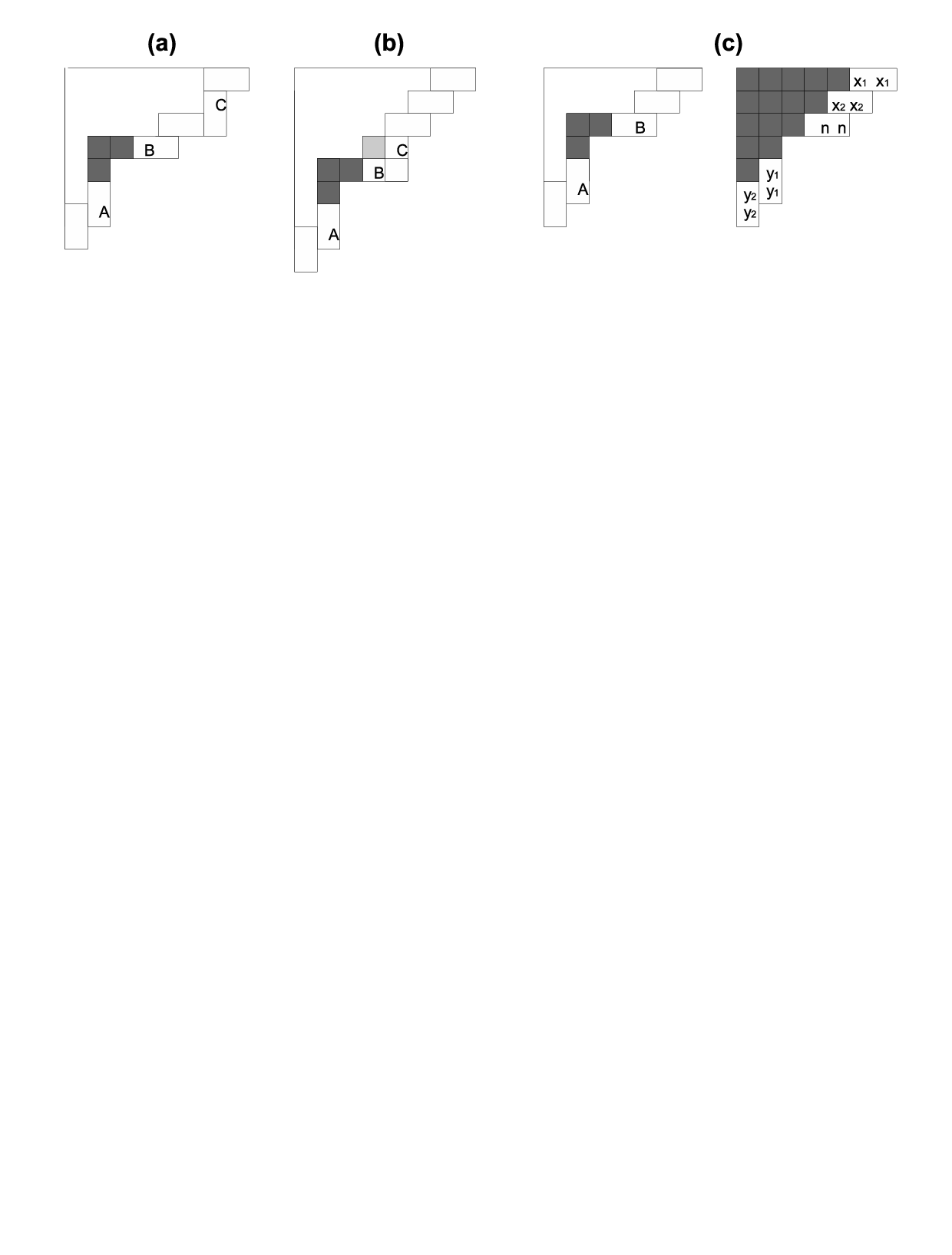}
\end{array}$
\end{center} \vspace{-5in}
\caption{Case 4. }\label{main.theorem.fig5}
\end{figure}

\vskip.1in \noindent {\bf Case 4.1:} We first assume that   there is another  domino corner $C$ of $T$.  W.L.O.G. we  assume that $C$  lies  in
$(T,B, \mathrm{ne})$ since the other case can be dealt with in the same manner after taking the transpose of $T$.
\vskip.1in \noindent {\bf Case 4.1.1:} Suppose that
$(\mathrm{T,C,sw})\cap(T,B,\mathrm{ne})$  contains a domino corner
as in Figure~\ref{main.theorem.fig5}($a$). Let $\sigma$ be a signed
permutation such that $P(\sigma)=T^{\uparrow C}$. Then the pairs of
domino corners $A$ and $C$ and similarly $B$ and $C$ satisfy  Case
2. Therefore we have $\beta ~\rp~\sigma~\rp~ \alpha$. The case when
$C$ lies in $(T,A,\mathrm{sw})$ follows similarly.

\vskip.1in \noindent {\bf Case 4.1.2:}  Now suppose that
$(\mathrm{T,C,sw})\cap(T,B,\mathrm{ne})$ is a staircase shape
$(s,s-1,\ldots,2,1)$ for some $s\geq 1$. Observe that the case
 $s>1$ is impossible since $B$ is a horizontal domino cell.
For  $s=1$ consider  Figure~\ref{main.theorem.fig5}($b$).
Let $\sigma$ be a sign permutation  such that $ T^{\uparrow
C}=P^r(\sigma)$. Then as the domino corners $B$ and $C$ satisfy Case
3.1.1, we have
 $\beta ~\rp~ \sigma$.  On the other hand
$(T,A,\mathrm{ne})\cap(T,C,\mathrm{sw})$ also contains a domino corner,
 therefore by Case 2 we have $\alpha ~\rp~ \sigma$. Hence
$\alpha ~\rp~ \beta$. The case when $C$ lies in $(T,A,\mathrm{sw})$
also  follows similarly.

\vskip.1in \noindent {\bf Case 4.2:} Now we suppose that there is no domino corner  of $T$ beyond $A$ and
$B$ as illustrated in Figure~\ref{main.theorem.fig5}($c$). One can
easily see that after reverse insertion all horizontal domino cells
and then vertical domino cells, only a staircase shape is left.
Therefore $n\leq r+1$ and $\eta(T^{\uparrow A})= \alpha_n<0$ and
$\eta(T^{\uparrow B})=\beta_n>0$. Moreover we have
$$
T^{\uparrow A\uparrow B}=T^{\uparrow B \uparrow A},~
\eta(T^{\uparrow A \uparrow B})=\beta_n ~\mbox{ and }~
\eta(T^{\uparrow B \uparrow A})=\alpha_n.$$

 Let $u$ be a signed word such that $P^r(u)=T^{\uparrow A\uparrow
B}=T^{\uparrow B\uparrow A}$. Clearly $P^r(u\alpha_n\beta_n)=T=
P^r(u\beta_n\alpha_n)$ and the size of $u$ is less than $r-1$.
Moreover
$$u\alpha_n\beta_n ~\rdb~u \beta_n\alpha_n.$$ On the
other hand $P^r(u\beta_n)=T^{\uparrow A}$ and
$P^r(u\alpha_n)=T^{\uparrow B}$ and by  induction we have
$\alpha_1\ldots\alpha_{n-1}~\rp~u\beta_n$ and
$\beta_1\ldots\beta_{n-1}~\rp~ u\alpha_n$. Hence
$\alpha=\alpha_1\ldots\alpha_{n-1}\alpha_n ~\rp~ u\beta_n\alpha_n
~\rdb~u\alpha_n\beta_n ~\rp~ \beta_1\ldots\beta_{n-1}\beta_n=\beta$
as desired.

\end{proof}

 \end{document}